%% file: BinJiaThesis.tex
\documentclass[a4paper, 12pt, twoside]{book}

\usepackage[utf8]{inputenc}
\usepackage{graphicx}
\usepackage{fancyhdr}
\usepackage[fleqn]{amsmath}
\usepackage{accents}
\usepackage{stfloats}
\usepackage{setspace}
\usepackage{amsthm}
\usepackage{lineno}

\usepackage[Lenny]{fncychap}   

\usepackage{mathtools}
\usepackage{mathrsfs}
\usepackage{bbm}
\usepackage{latexsym,amscd,amsfonts,epsfig,graphicx}

\usepackage{ifpdf, color}%

\usepackage{stmaryrd}
\usepackage{pifont}
\usepackage{mathbbold}
\usepackage{amsmath, mathtools}
\usepackage{mathrsfs}
\usepackage{bbm}
\usepackage{amsfonts}
\usepackage{latexsym,amssymb,amsmath,amscd,amsfonts,epsfig,graphicx}

 \allowdisplaybreaks

\usepackage{subfigure}
\usepackage{palatino}
\usepackage{xr}
\usepackage{dsfont}

\usepackage{charter}

\usepackage{multirow}

\usepackage[colorlinks=true, linkcolor=black, citecolor=black,urlcolor=blue,
pdftitle={A study of link graphs},
pdfauthor={Bin JIA},
pdfsubject={Thesis},
pdfkeywords={graph, link graph}]{hyperref}

\long\def\delete#1{}

\usepackage[lined,commentsnumbered,algosection,ruled,vlined]{algorithm2e}

\newtheorem{thm}{Theorem}[section]
\newtheorem{lem}[thm]{Lemma}

\newtheorem{coro}[thm]{Corollary}
\newtheorem{prop}[thm]{Property}
\newtheorem{propo}[thm]{Proposition}
\newtheorem{obs}[thm]{Observation}

\def\pf{\noindent{\it Proof.\;\;\:}}
\def\qed{\nopagebreak\hfill{\rule{4pt}{7pt}}
\medbreak}

\DeclareMathOperator{\diam}{diam}
\DeclareMathOperator{\dege}{d}

\DeclareMathOperator{\radi}{radi}
\DeclareMathOperator{\dist}{dist}
\DeclareMathOperator{\ecc}{ecc}
\DeclareMathOperator{\tdi}{tdi}
\DeclareMathOperator{\tw}{tw}
\DeclareMathOperator{\girth}{girth}
\theoremstyle{definition}
\newtheorem{defi}[thm]{Definition}
\newtheorem{exa}[thm]{Example}

\usepackage{bbding}
\usepackage{pifont}
\usepackage{txfonts}
\usepackage{stmaryrd}

\usepackage{amssymb}

\begin{document}

\frontmatter



\include{Title}

\include{Abstract}

\tableofcontents

\listoffigures

%
\mainmatter

\include{Introduction}
\include{Chapter11}

\include{Chapter21}

\include{Sec31}

\include{Sec32}
\include{Sec33}

\include{Sec35}

\include{Sec36}

\include{41Introduction}

\include{42Terminology}
\include{43Examples}

\include{44IncidentL-links}

\include{45partitionedGraphs}

\include{46BoundingSize}

\include{47BoundingMaximumDegree}
\include{48PathGraphs}
\cleardoublepage

\bibliographystyle{plain}
\bibliography{Xbib}   
\nocite{}
\clearpage
\newpage




\end{document}

%% file: Title.tex
\begin{titlepage}

\hfill\\
\hfill\\
\vspace{1cm}

\begin{center}
\textbf{\LARGE A study of link graphs}

\vspace{3cm}

\textbf{\normalsize Contents are presented at\\
Discrete Maths Research Group Seminar\\
Monash University, 30 May 2014
}


\vspace{3cm}

\textbf{\large Bin Jia}

\medskip

\vspace{40pt}

Department of Mathematics and Statistics\medskip

The University of Melbourne\medskip

Australia\medskip


\end{center}
\end{titlepage}


       \newpage

   ~
        \vspace{15cm}

\begin{center}
 \copyright ~   Copyright by Bin Jia, 2014. ~~~~~~~~~~~~
\end{center}

%% file: Abstract.tex
\addcontentsline{toc}{chapter}{Abstract}
\vspace{4cm}
\begin{center}
\textbf{\large Abstract}
\end{center}


This thesis contributes in various aspects to the characterisation and determination problems for incidence patterns proposed by Gr{\"u}nbaum (1969).
Specifically, we introduce and study a new kind of incidence pattern, called \emph{$\ell$-link graphs},
which generalises the notions of line graphs and path graphs. An \emph{$\ell$-link}
is a walk of length
$\ell \geqslant 0$ in that graph such that
consecutive edges are different. We identify an $\ell$-link with its reverse sequence. For example, a $0$-link is a vertex. And a $1$-link is an edge. Further, an $\ell$-path is an $\ell$-link without repeated vertices. The
\emph{$\ell$-link graph $\mathbb{L}_\ell(G)$}
of a graph $G$ is the graph with vertices the $\ell$-links of $G$, such that two vertices are adjacent if the union of their corresponding $\ell$-links is an $(\ell + 1)$-link; Or equivalently, one corresponding $\ell$-link can be shunted to the other in one step. The definition here is for simple graphs, but will be extended to graphs with parallel edges.

We reveal a recursive structure for $\ell$-link graphs, which allows us to bound the chromatic number of $\mathbb{L}_\ell(G)$ in terms of $\ell$ and the chromatic number or edge chromatic number of $G$. As a corollary, $\mathbb{L}_\ell(G)$ is $3$-colourable for each finite graph $G$ and large enough $\ell$. By investigating the shunting of $\ell$-links in $G$, we show that the Hadwiger number of a nonempty $\mathbb{L}_\ell(G)$ is at least that of $G$. Hadwiger's conjecture states that the Hadwiger number of a graph is at least the chromatic number of that graph. The conjecture has been proved by Reed and Seymour (2004) for line graphs, and hence $1$-link graphs. We prove the conjecture for a wide class of $\ell$-link graphs.

An \emph{$\ell$-root} of a graph $H$ is a graph $G$ such that $H \cong \mathbb{L}_\ell(G)$. For instance, $K_3$ and $K_{1, 3}$ are $1$-roots of $K_3$.
We show that every $\ell$-root of a finite graph is a certain combination of a finite minimal (up to the subgraph relation) $\ell$-root and trees of bounded diameter. This transfers the study of $\ell$-roots into that of finite minimal $\ell$-roots. As a generalisation of Whitney's theorem (1932), we bound from above the number, size, order and maximum degree of minimal $\ell$-roots of finite graphs.
This implies that the $\ell$-roots of a finite graph are better-quasi-ordered by the induced subgraph relation. This work forms the basis for solving the recognition and determination problems for $\ell$-link graphs in our future papers. Similar results are obtained for path graphs (Broersma and Hoede, 1989). $G$ is an \emph{$\ell$-path root} of a graph $H$ if $H$ is isomorphic to the $\ell$-path graph of $G$. We bound from above the number, size and order of minimal $\ell$-path roots of a finite graph. Further, we show that every sequence of $\ell$-path roots of a finite graph is better-quasi-ordered by the subgraph relation, and by the induced subgraph relation if these roots have bounded multiplicity.

%% file: Chapter11.tex
\chapter{Introduction}
\section{Incidence patterns}
Introduced by Gr{\"u}nbaum \cite{Grunbaum1969}, an \emph{incidence pattern} is a function that maps given graphs or similar objects to graphs. Two general questions associated with a given incidence pattern were proposed by Gr{\"u}nbaum as characterising all graphs that can be constructed from this pattern, and determining the original object for each of these graphs. This thesis contributes in various aspects to the two questions for a certain incidence pattern, called the \emph{$\ell$-link graph construction} \cite{JiaWood2013}.

%% file: Chapter21.tex
\chapter{Preliminaries}

\section{Notation, definitions and preliminaries}

%% file: Sec31.tex
\section{Introduction and main results}\label{sec_intro}
We introduce a new family of graphs, called \emph{$\ell$-link graphs}, which generalises the notions of line graphs and path graphs. Such a graph is constructed from a certain kind of walk of length $\ell \geqslant 0$ in a given graph $G$. To ensure that the constructed graph is undirected, $G$ is  undirected, and we identify a walk with its reverse sequence. To avoid loops, $G$ is loopless, and the consecutive edges in each walk are different. Such a walk is called an \emph{$\ell$-link}. For example, a $0$-link is a vertex, a $1$-link is an edge, and a $2$-link consists of two edges with an end vertex in common. An \emph{$\ell$-path} is an $\ell$-link without repeated vertices. We use $\mathscr{L}_\ell(G)$ and $\mathscr{P}_\ell(G)$ to denote the sets of $\ell$-links and $\ell$-paths of $G$ respectively. There have been a number of families of graphs constructed from $\ell$-links. As one of the most commonly studied graphs, the \emph{line graph $\mathbb{L}(G)$}, introduced by Whitney \cite{Whitney1932}, is the simple graph with vertex set $E(G)$, in which two vertices are adjacent if their corresponding edges are incident to a common vertex. More generally, the \emph{$\ell$-path graph $\mathbb{P}_{\ell}(G)$} is the simple graph with vertex set $\mathscr{P}_{\ell}(G)$, where two vertices are adjacent if the union of their corresponding $\ell$-paths forms a path or a cycle of length $\ell + 1$. Note that $\mathbb{P}_{\ell}(G)$ is the $\mathbb{P}_{\ell + 1}$-graph of $G$ introduced by Broersma and Hoede \cite{BH1989}. Inspired by these graphs, we define the \emph{$\ell$-link graph $\mathbb{L}_\ell(G)$} of $G$ to be the graph with vertex set $\mathscr{L}_\ell(G)$, in which two vertices are joined by $\mu \geqslant 0$ edges in $\mathbb{L}_\ell(G)$ if they correspond to two subsequences of each of $\mu$ $(\ell + 1)$-links of $G$. More strict definitions can be found in Section \ref{sec_defiTermi}, together with some other related graphs.

This paper studies the structure, colouring and minors of $\ell$-link graphs including a proof of Hadwiger's conjecture for a wide class of $\ell$-link graphs. By default $\ell \geqslant 0$ is an integer. And all graphs are finite, undirected and loopless. Parallel edges are admitted unless we specify the graph to be \emph{simple}.

\subsection{Graph colouring} Let $t \geqslant 0$ be an integer. A \emph{$t$-colouring} of $G$ is a map $\lambda: V(G) \rightarrow [t] := \{1, 2, \ldots, t\}$ such that $\lambda(u) \neq \lambda(v)$ whenever $u, v \in V(G)$ are adjacent in $G$. A graph with a $t$-colouring is \emph{$t$-colourable}. The \emph{chromatic number $\chi(G)$} is the minimum $t$ such that $G$ is $t$-colourable. Similarly, an \emph{$t$-edge-colouring} of $G$ is a map $\lambda: E(G) \rightarrow [t]$ such that $\lambda(e) \neq \lambda(f)$ whenever $e, f \in E(G)$ are incident to a common vertex in $G$. The \emph{edge-chromatic number $\chi'(G)$} of $G$ is the minimum $t$ such that $G$ admits a $t$-edge-colouring. Let $\chi_\ell(G) := \chi(\mathbb{L}_\ell(G))$, and \emph{$\Delta(G)$} be the maximum degree of $G$. By \cite[Proposition 5.2.2]{Diestel2010}, $\chi_0(G) = \chi(G) \leqslant \Delta(G) + 1$. Shannon \cite{Shannon1949} proved that $\chi_1(G) = \chi'(G) \leqslant \frac{3}{2}\Delta(G)$. We prove a recursive structure for $\ell$-link graphs which leads to the following upper bounds for $\chi_\ell(G)$:

\begin{thm}\label{thm_chrl}
Let $G$ be a graph, $\chi := \chi(G)$, $\chi' := \chi'(G)$, and $\Delta := \Delta(G)$.
\begin{itemize}
\item[\bf (1)] If $\ell \geqslant 0$ is even, then $\chi_\ell(G) \leqslant \min\{\chi, \lfloor(\frac{2}{3})^{\ell/2}(\chi - 3)\rfloor + 3\}$.
\item[\bf (2)] If $\ell \geqslant 1$ is odd, then $\chi_\ell(G) \leqslant \min\{\chi', \lfloor(\frac{2}{3})^{\frac{\ell - 1}{2}}(\chi' - 3)\rfloor + 3\}$.
\item[\bf (3)] If $\ell \neq 1$, then $\chi_\ell(G) \leqslant \Delta + 1$. 
\item[\bf (4)] If $\ell \geqslant 2$, then $\chi_\ell(G) \leqslant \chi_{\ell - 2}(G)$. 
\end{itemize}
\end{thm}

Theorem \ref{thm_chrl} implies that $\mathbb{L}_\ell(G)$ is $3$-colourable for large enough $\ell$.

\begin{coro}\label{coro_chrl} For each graph $G$,
$\mathbb{L}_\ell(G)$ is $3$-colourable in the following cases:
\begin{itemize}
\item[{\bf (1)}] $\ell \geqslant 0$ is even, and either $\chi(G) \leqslant 3$ or $\ell > 2\log_{1.5}(\chi(G) - 3)$.
\item[{\bf (2)}] $\ell \geqslant 1$ is odd, and either $\chi'(G) \leqslant 3$ or $\ell > 2\log_{1.5}(\chi'(G) - 3) + 1$.
\end{itemize}
\end{coro}

As explained in Section \ref{sec_defiTermi}, this corollary is related to and implies a result by Kawai and Shibata \cite{KawaiShibata2002}.

\subsection{Graph minors} By \emph{contracting} an edge we mean identifying its end vertices and deleting possible resulting loops. A graph $H$ is a \emph{minor} of $G$ if $H$ can be obtained from a subgraph of $G$ by contracting edges. An \emph{$H$-minor} is a minor of $G$ that is isomorphic to $H$. The \emph{Hadwiger number $\eta(G)$} of $G$ is the maximum integer $t$ such that $G$ contains a $K_t$-minor. Denote by \emph{$\delta(G)$} the minimum degree of $G$. The \emph{degeneracy $\dege(G)$} of $G$ is the maximum $\delta(H)$ over the subgraphs $H$ of $G$. We prove the following:
\begin{thm}\label{thm_hgeqg}
Let $\ell \geqslant 1$, and $G$ be a graph such that $\mathbb{L}_\ell(G)$ contains at least one edge. Then $\eta(\mathbb{L}_\ell(G)) \geqslant \max\{\eta(G), \dege(G)\}$.
\end{thm}

By definition $\mathbb{L}(G)$ is the underlying simple graph of $\mathbb{L}_1(G)$. And $\mathbb{L}_\ell(G) = \mathbb{P}_\ell(G)$ if $girth(G) > \{\ell, 2\}$.
Thus Theorem \ref{thm_hgeqg} can be applied to path graphs.
\begin{coro}\label{coro_hgeqg}
Let $\ell \geqslant 1$, and $G$ be a graph of girth at least $\ell + 1$ such that $\mathbb{P}_\ell(G)$ contains at least one edge. Then
$\eta(\mathbb{P}_\ell(G)) \geqslant \max\{\eta(G), \dege(G)\}$.
\end{coro}

As a far-reaching generalisation of the four-colour theorem, in 1943, Hugo Hadwiger \cite{Hadwiger1943} conjectured the following:

\medskip
\noindent {\bf Hadwiger's conjecture:} $\eta(G) \geqslant \chi(G)$ for every graph $G$.
\medskip

Hadwiger's conjecture was proved by Robertson, Seymour and Thomas  \cite{RST1993} for $\chi(G) \leqslant 6$. The conjecture for line graphs, or equivalently for $1$-link graphs, was proved by Reed and Seymour \cite{Reed2004}. We prove the following:
\begin{thm}\label{thm_hcllink}
Hadwiger's conjecture is true for $\mathbb{L}_\ell(G)$ in the following cases:
\begin{itemize}
\item[{\bf (1)}] $\ell \geqslant 1$ and $G$ is biconnected. 
\item[{\bf (2)}] $\ell \geqslant 2$ is an even integer. 
\item[{\bf (3)}] $\dege(G) \geqslant 3$ and $\ell > 2\log_{1.5}\frac{\Delta(G) - 2}{\dege(G) - 2} + 3$. 
\item[{\bf (4)}] $\Delta(G) \geqslant 3$ and $\ell > 2\log_{1.5}(\Delta(G) - 2) - 3.83$. 
\item[{\bf (5)}] $\Delta(G) \leqslant 5$.
\end{itemize}
\end{thm}

The corresponding results for path graphs are listed below:
\begin{coro}
Let $G$ be a graph of girth at least $\ell + 1$. Then Hadwiger's conjecture holds for $\mathbb{P}_{\ell}(G)$ in the cases of Theorem \ref{thm_hcllink} (1) -- (5).
\end{coro}

%% file: Sec32.tex
\section{Definitions and terminology}\label{sec_defiTermi}
We now give some formal definitions. A graph $G$ is \emph{null} if $V(G) = \emptyset$, and \emph{nonnull} otherwise. A nonnull graph $G$ is \emph{empty} if $E(G) = \emptyset$, and nonempty otherwise. A \emph{unit} is a vertex or an edge. The subgraph of $G$ induced by $V \subseteq V(G)$ is the maximal subgraph of $G$ with vertex set $V$. And in this case, the subgraph is called an \emph{induced subgraph} of $G$. For $\emptyset \neq E \subseteq E(G)$, the subgraph of $G$ induced by $E \cup V$ is the minimal subgraph of $G$ with edge set $E$, and vertex set including $V$.

For more accurate analysis, we need to define \emph{$\ell$-arcs}. An \emph{$\ell$-arc} (or \emph{$*$-arc} if we ignore the length) of $G$ is an alternating sequence $\vec{L} := (v_0, e_1, \ldots, e_{\ell}, v_\ell)$ of units of $G$ such that the end vertices of $e_i \in E(G)$ are $v_{i - 1}$ and $v_i$ for $i \in [\ell]$, and that $e_i \neq e_{i + 1}$ for $i \in [\ell - 1]$. The \emph{direction} of $\vec{L}$ is its vertex sequence \emph{$(v_0, v_1, \ldots, v_\ell)$}. In algebraic graph theory, $\ell$-arcs in simple graphs have been widely studied \cite{Tutte1947,Tutte1959,Weiss1981,Biggs1993}. Note that $\vec{L}$ and its \emph{reverse $-\vec{L} := (v_\ell, e_{\ell}, \ldots, e_1, v_0)$} are different unless $\ell = 0$. The $\ell$-link (or \emph{$*$-link} if the length is ignored) $L := [v_0, e_1, \ldots, e_{\ell}, v_\ell]$ is obtained by taking $\vec{L}$ and $-\vec{L}$ as a single object. For $0 \leqslant i \leqslant j \leqslant \ell$, the $(j - i)$-arc \emph{$\vec{L}(i, j) := (v_i, e_{i + 1},\ldots, e_{j}, v_j)$} and the $(j - i)$-link \emph{$\vec{L}[i, j] := [v_i, e_{i + 1}, \ldots, e_{j}, v_j]$} are called \emph{segments} of $\vec{L}$ and $L$ respectively. We may write $\vec{L}(j, i) := -\vec{L}(i, j)$, and $\vec{L}[j, i] := \vec{L}[i, j]$. These segments are called \emph{middle segments} if $i + j = \ell$. $L$ is called an \emph{$\ell$-cycle} if $\ell \geqslant 2$, $v_0 = v_\ell$ and $\vec{L}[0, \ell - 1]$ is an $(\ell - 1)$-path. Denote by $\vec{\mathscr{L}}_{\ell}(G)$ and $\mathscr{C}_\ell(G)$ the sets of $\ell$-arcs and $\ell$-cycles of $G$ respectively. Usually, $\vec{e}_i := (v_{i - 1}, e_i, v_i)$ is called an \emph{arc} for short. In particular, $v_0$, $v_\ell$, $e_1$, $e_\ell$, $\vec{e}_1$ and $\vec{e}_\ell$ are called the \emph{tail vertex}, \emph{head vertex}, \emph{tail edge}, \emph{head edge}, \emph{tail arc}, and \emph{head arc} of $\vec{L}$ respectively.

\begin{figure}[!h]
\centering
\includegraphics{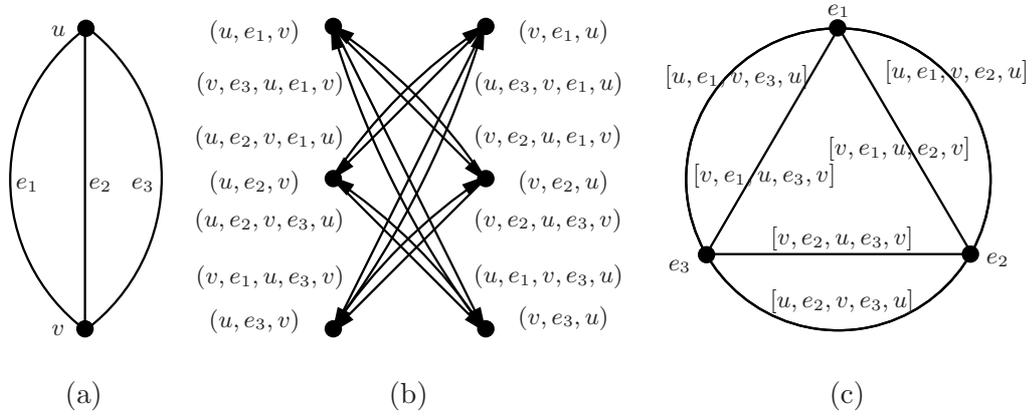}
\caption{(a) $D_3$ $\quad$ (b) $\mathbb{A}_1(D_3)$ $\quad$ (c) $\mathbb{L}_1(D_3)$}
\label{F:D3A1L1}
\end{figure}

Godsil and Royle \cite{GR2001} defined the \emph{$\ell$-arc graph $\mathbb{A}_{\ell}(G)$} to be the digraph with vertex set $\vec{\mathscr{L}}_\ell(G)$, such that there is an arc, labeled by $\vec{Q}$, from $\vec{Q}(0, \ell)$ to $\vec{Q}(1, \ell + 1)$ in $\mathbb{A}_{\ell}(G)$ for every $\vec{Q} \in \vec{L}_{\ell + 1}(G)$. The \emph{$t$-dipole graph $D_t$} is the graph consists of two vertices and $t \geqslant 1$ edges between them. (See Figure \ref{F:D3A1L1}(a) for $D_3$, and Figure \ref{F:D3A1L1}(b) the $1$-arc graph of $D_3$.)
The \emph{$\ell^{th}$ iterated line digraph $\mathbb{A}^\ell(G)$} is $\mathbb{A}_1(G)$ if $\ell = 1$, and $\mathbb{A}_1(\mathbb{A}^{\ell - 1}(G))$ if $\ell \geqslant 2$ (see  \cite{BW1978}). Examples of undirected graphs constructed from $\ell$-arcs can be found in \cite{JLW2010,Jia2011}.

\begin{figure}[!h]
\centering
\includegraphics{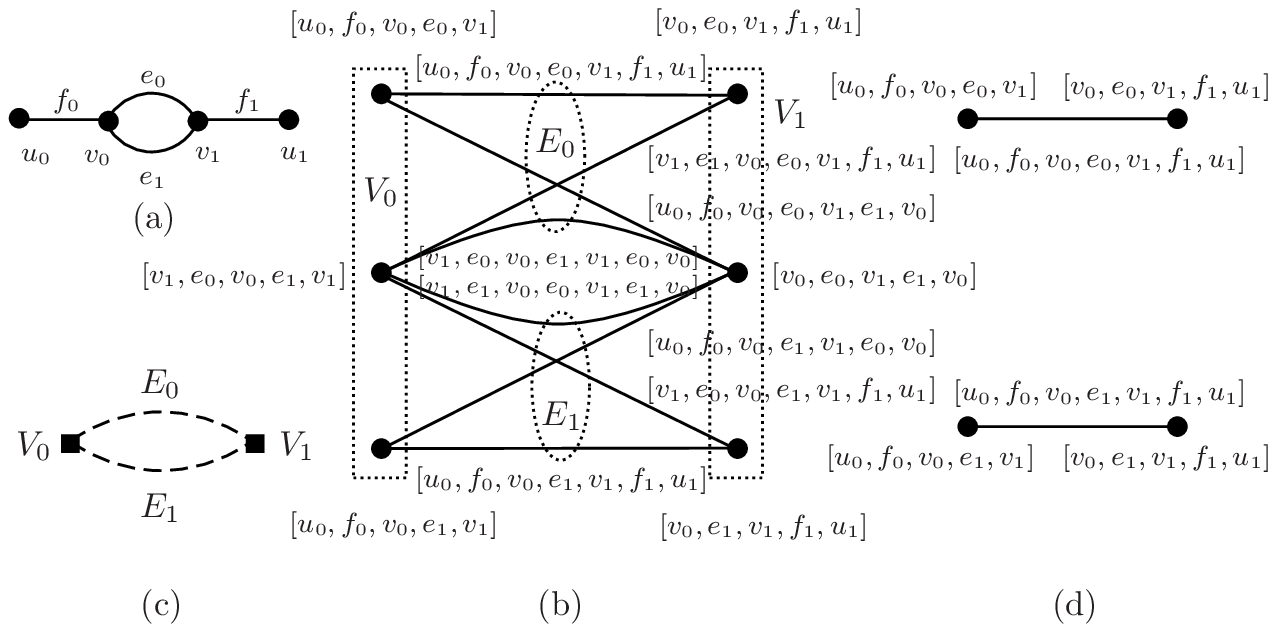}
\caption{(a) $G$ $\quad$ (b) $H := \mathbb{L}_2(G)$ $\quad$ (c) $H_{(\mathcal{V}, \mathcal{E})}$ $\quad$ (d) $\mathbb{P}_2(G)$}
\label{F:GL2P2Quotient}
\end{figure}

\emph{Shunting} of $\ell$-arcs was introduced by Tutte \cite{Tutte1966}. We extend this motion to $\ell$-links. For $\ell, s \geqslant 0$, and $\vec{Q} \in \vec{\mathscr{L}}_{\ell + s}(G)$, let $\vec{L}_i := \vec{Q}(i, \ell + i)$ for $i \in [0, s]$, and $\vec{Q}_i := \vec{L}(i - 1, \ell + i)$ for $i \in [s]$. Let \emph{$Q^{[\ell]} := [L_0, Q_1, L_1, \ldots, L_{s - 1}, Q_{s}, L_{s}]$}. We say $L_0$ can be \emph{shunted} to $L_s$ through $\vec{Q}$ or $Q$. \emph{$Q^{\{\ell\}} := \{L_0, L_1, \ldots, L_s\}$} is the set of \emph{images} during this shunting. For $L, R \in \mathscr{L}_\ell(G)$, we say $L$ can be \emph{shunted} to $R$ if there are $\ell$-links $L = L_0, L_1, \ldots, L_s = R$ such that $L_{i - 1}$ can be shunted to $L_{i}$ through some $*$-arc $\vec{Q}_i$ for $i \in [s]$. In Figure \ref{F:GL2P2Quotient}, $[u_0, f_0, v_0, e_0, v_1]$ can be shunted to $[v_1, e_0, v_0, e_1, v_1]$ through $(u_0, f_0, v_0, e_0, v_1, f_1, u_1)$ and $(u_1, f_1, v_1, e_0, v_0, e_1, v_1)$.

For $L, R \in \mathscr{L}_\ell(G)$ and $\mathscr{Q} \subseteq \mathscr{L}_{\ell + 1}(G)$, denote by \emph{$\mathscr{Q}(L, R)$} the set of $Q \in \mathscr{Q}$ such that $L$ can be shunted to $R$ through $Q$. We show in Section \ref{sec_genStruc} that $|\mathscr{Q}(L, R)|$ is $0$ or $1$ if $G$ is simple, and can be up to $2$ if $\ell \geqslant 1$ and $G$ contains parallel edges. A more formal definition of $\ell$-link graphs is given below:
\begin{defi}\label{defi_llinkgraph}
Let $\mathscr{L} \subseteq \mathscr{L}_\ell(G)$, and $\mathscr {Q} \subseteq \mathscr{L}_{\ell + 1}(G)$. The \emph{partial $\ell$-link graph $\mathbb{L}(G, \mathscr{L}, \mathscr{Q})$}  of $G$, with respect to $\mathscr{L}$ and $\mathscr{Q}$, is the graph with vertex set $\mathscr{L}$, such that $L, R \in \mathscr{L}$ are joined by exactly $|\mathscr{Q}(L, R)|$ edges. In particular, $\mathbb{L}_\ell(G) = \mathbb{L}(G, \mathscr{L}_\ell(G), \mathscr{L}_{\ell + 1}(G))$ is the $\ell$-link graph of $G$.
\end{defi}
\noindent {\bf Remark.} We assign exclusively to each edge of $\mathbb{L}_\ell(G)$ between $L, R \in \mathscr{L}_\ell(G)$ a $Q \in \mathscr{L}_{\ell + 1}(G)$ such that $L$ can be shunted to $R$ through $Q$, and refer to this edge simply as $Q$. In this sense, $Q^{[\ell]} := [L, Q, R]$ is a $1$-link of $\mathbb{L}_{\ell}(G)$.

For example, the $1$-link graph of $D_3$ can be seen in Figure \ref{F:D3A1L1}(c). A $2$-link graph is given in Figure \ref{F:GL2P2Quotient}(b), and a $2$-path graph is depicted in Figure \ref{F:GL2P2Quotient}(d).

Reed and Seymour \cite{Reed2004} pointed out that proving Hadwiger's conjecture for line graphs of multigraphs is more difficult than for that of simple graphs. This motivates us to work on the $\ell$-link graphs of multigraphs. Diestel \cite[page 28]{Diestel2010} explained that, in some situations, it is more natural to develop graph theory for multigraphs. We allow parallel edges in $\ell$-link graphs in order to give a characterisation for $\mathbb{L}_\ell(G)$ in our companion papers \cite{JiaWood2013lleq3,JiaWood2013lgeq4} regardless of whether $G$ is simple. The observation below follows from the definitions:

\begin{obs}
$\mathbb{L}_0(G) = G$, $\mathbb{P}_1(G) = \mathbb{L}(G)$, and $\mathbb{P}_{\ell}(G)$ is the underlying simple graph of $\mathbb{L}_\ell(G)$ for $\ell \in \{0, 1\}$. For $\ell \geqslant 2$, $\mathbb{P}_{\ell}(G) = \mathbb{L}(G, \mathscr{P}_{\ell}(G), \mathscr{P}_{\ell + 1}(G)\\ \cup \mathscr{C}_{\ell + 1}(G))$ is an induced subgraph of $\mathbb{L}_{\ell}(G)$. If $G$ is simple, then $\mathbb{P}_{\ell}(G) = \mathbb{L}_\ell(G)$ for $\ell \in \{0, 1, 2\}$. Further, $\mathbb{P}_{\ell}(G) = \mathbb{L}_\ell(G)$ if $girth(G) > \max\{\ell, 2\}$.
\end{obs}


Let $\vec{Q} \in \vec{\mathscr{L}}_{\ell + s}(G)$, and $[L_0, Q_1, L_1, \ldots, L_{s - 1}, Q_{s}, L_{s}] := Q^{[\ell]}$. From Definition \ref{defi_llinkgraph}, for $i \in [s]$, $Q_i$ is an edge of $H := \mathbb{L}_\ell(G)$ between $L_{i - 1}, L_i \in V(H)$. So $Q^{[\ell]}$ is an $s$-link of $H$. In Figure \ref{F:GL2P2Quotient}(b), $[u_0, f_0, v_0, e_0, v_1, e_1, v_0, e_0, v_1]^{[2]} = [[u_0, f_0, v_0, e_0, v_1], [u_0, f_0, v_0, e_0, v_1, e_1, v_0], [v_0, e_0, v_1, e_1, v_0], [v_0, e_0, v_1, e_1, v_0, e_0, v_1]$, $[v_1, e_1, v_0, e_0, v_1]]$ is a $2$-path of $H$.

We say $H$ is \emph{homomorphic} to $G$, written $H \rightarrow G$, if there is an injection $\alpha: V(H) \cup E(H) \rightarrow V(G) \cup E(G)$ such that for $w \in V(H)$, $f \in E(H)$ and $[u, e, v] \in \mathscr{L}_1(H)$, their images $w^\alpha \in V(G)$, $f^\alpha \in E(G)$ and $[u^\alpha, e^\alpha, v^\alpha] \in \mathscr{L}_1(G)$. In this case,  $\alpha$ is called a \emph{homomorphism} from $H$ to $G$. The definition here is a generalisation of the one for simple graphs by Godsil and Royle \cite[Page 6]{GR2001}. A bijective homomorphism is an \emph{isomorphism}. By Hell and Ne{\v{s}}et{\v{r}}il \cite{HellNesetril2004}, $\chi(H) \leqslant \chi(G)$ if $H \rightarrow G$. For instance, $\vec{L} \mapsto L$ for $\vec{L} \in \vec{\mathscr{L}}_{\ell}(G) \cup \vec{\mathscr{L}}_{\ell + 1}(G)$ can be seen as a homomorphism from $\mathbb{A}_\ell(G)$ to $\mathbb{L}_\ell(G)$. By Bang-Jensen and Gutin \cite{BangGutin2009}, $\mathbb{A}_\ell(G) \cong \mathbb{A}^\ell(G)$. So $\chi(\mathbb{A}^\ell(G)) = \chi(\mathbb{A}_\ell(G)) \leqslant \chi(\mathbb{L}_\ell(G)) \leqslant \chi_\ell(G)$. We emphasize that $\chi(\mathbb{A}^\ell(G))$ might be much less than $\chi_\ell(G)$. For example, as depicted in Figure \ref{F:D3A1L1}, when $t \geqslant 3$, $\chi(\mathbb{A}^\ell(D_t)) = 2 < t = \chi_\ell(D_t)$. Kawai and Shibata proved that $\mathbb{A}^\ell(G)$ is $3$-colourable for large enough $\ell$. By the analysis above, Corollary \ref{coro_chrl} implies this result.


A graph homomorphism from $H$ is usually represented by a vertex partition $\mathcal {V}$ and an edge partition $\mathcal {E}$ of $H$ such that: {\bf (a)} each part of $\mathcal {V}$ is an independent set of $H$, and {\bf (b)} each part of $\mathcal {E}$ is incident to exactly two parts of $\mathcal {V}$. In this situation, for different $U, V \in \mathcal {V}$, define \emph{$\mu(U, V)$} to be the number of parts of $\mathcal {E}$ incident to both $U$ and $V$. The \emph{quotient graph $H_{(\mathcal {V}, \mathcal {E})}$} of $H$ is defined to be the graph with vertex set $\mathcal {V}$, and for every pair of different $U, V \in \mathcal {V}$, there are exactly $\mu(U, V)$ edges between them. To avoid ambiguity, for $V \in \mathcal {V}$ and $E \in \mathcal {E}$, we use $V_{\mathcal {V}}$ and $E_{\mathcal {E}}$ to denote the corresponding vertex and edge of $H_{(\mathcal {V}, \mathcal {E})}$, which defines a graph homomorphism from $H$ to $H_{(\mathcal {V}, \mathcal {E})}$. Sometimes, we only need the underlying simple graph \emph{$H_{\mathcal {V}}$} of $H_{(\mathcal {V}, \mathcal {E})}$.

For $\ell \geqslant 2$, there is a natural partition in an $\ell$-link graph. For each $R \in \mathscr{L}_{\ell - 2}(G)$, let $\mathscr{L}_\ell(R)$ be the set of $\ell$-links of $G$ with middle segment $R$. Clearly, $\mathcal {V}_{\ell}(G) := \{\mathscr{L}_{\ell}(R) \neq \emptyset | R \in \mathscr{L}_{\ell - 2}(G)\}$ is a vertex partition of $\mathbb{L}_\ell(G)$. And $\mathcal {E}_\ell(G) := \{\mathscr{L}_{\ell + 1}(P) \neq \emptyset | P \in \mathscr{L}_{\ell - 1}(G)\}$ is an edge partition of $\mathbb{L}_\ell(G)$. Consider the $2$-link graph $H$ in Figure \ref{F:GL2P2Quotient}(b). The vertex and edge partitions of $H$ are indicated by the dotted rectangles and ellipses respectively. The corresponding quotient graph is given in Figure \ref{F:GL2P2Quotient}(c).

Special partitions are required to describe the structure of $\ell$-link graphs. Let $H$ be a graph admitting partitions $\mathcal {V}$ of $V(H)$ and $\mathcal {E}$ of $E(H)$ that satisfy {\bf (a)} and {\bf (b)} above. $(\mathcal {V}, \mathcal {E})$ is called an \emph{almost standard partition} of $H$ if further:

{\bf (c)} each part of $\mathcal {E}$ induces a complete bipartite subgraph of $H$,

{\bf (d)} each vertex of $H$ is incident to at most two parts of $\mathcal {E}$,

{\bf (e)} for each $V \in \mathcal{V}$, and different $E, F \in \mathcal {E}$, $V$ contains at most one vertex incident to both $E$ and $F$.

We use the term `almost standard partition' because the notion of `standard partition' is used in our companion paper \cite{JiaWood2013sdg3}.





%% file: Sec33.tex
\section{General structure of $\ell$-link graphs}\label{sec_genStruc}
We begin by determining some basic properties of $\ell$-link graphs, including their multiplicity and connectedness. The work in this section forms the basis for our main results on colouring and minors of $\ell$-link graphs.


Let us first fix some concepts by two observations.

\begin{obs}\label{obs_llinkregu}
The number of edges of $\mathbb{L}_{\ell}(G)$ is equal to the number of vertices of $\mathbb{L}_{\ell + 1}(G)$. In particular, if $G$ is $r$-regular for some $r \geqslant 2$, then this number is $|E(G)|(r - 1)^{\ell}$. If further $\ell \geqslant 1$, then $\mathbb{L}_{\ell}(G)$ is $2(r - 1)$-regular.
\end{obs}
\pf Let $G$ be $r$-regular, $n := |V(G)|$ and $m := |E(G)|$. We prove that $|\mathscr{L}_{\ell + 1}(G)| = m(r - 1)^{\ell}$ by induction on $\ell$. It is trivial for $\ell = 0$. For $\ell = 1$, $|\mathscr{L}_2(v)| = \binom{r}{2}$, and hence $|\mathscr{L}_{2}(G)| = \binom{r}{2}n = m(r - 1)$. Inductively assume $|\mathscr{L}_{\ell - 1}(G)| = m(r - 1)^{\ell - 2}$ for some $\ell \geqslant 2$. For each $R \in \mathscr{L}_{\ell - 1}(G)$, we have $|\mathscr{L}_{\ell + 1}(R)| = (r - 1)^2$ since $r \geqslant 2$. Thus $|\mathscr{L}_{\ell + 1}(G)| = |\mathscr{L}_{\ell - 1}(G)|(r - 1)^2 = m(r - 1)^{\ell}$ as desired. The other assertions follow from the definitions.
\qed

\begin{obs}\label{obs_llinkbipar}
Let $n, m \geqslant 2$. If $\ell \geqslant 1$ is odd, then $\mathbb{L}_\ell(K_{n, m})$ is $(n + m - 2)$-regular with order $nm[(n - 1)(m - 1)]^{\frac{\ell - 1}{2}}$. If $\ell \geqslant 2$ is even, then $\mathbb{L}_\ell(K_{n, m})$ has average degree $\frac{4(n - 1)(m - 1)}{n + m - 2}$, and order $\frac{1}{2}nm(n + m - 2)[(n - 1)(m - 1)]^{\frac{\ell}{2} - 1}$.
\end{obs}
\pf
Let $\ell \geqslant 1$ be odd, and $L$ be an $\ell$-link of $K_{n, m}$ with middle edge incident to a vertex $u$ of degree $n$ in $K_{n, m}$. It is not difficult to see that $L$ can be shunted in one step to $n - 1$ $\ell$-links whose middle edge is incident to $u$. By symmetry, each vertex of $\mathbb{L}_\ell(K_{n, m})$ is incident to $(n - 1) + (m - 1) = n + m - 2$ edges. Now we prove $|\mathscr{L}_\ell(K_{n, m})| = nm[(n - 1)(m - 1)]^{\frac{\ell - 1}{2}}$ by induction on $\ell$. Clearly, $|\mathscr{L}_1(K_{n, m})| = |E(K_{n, m})| = nm$. Inductively assume $|\mathscr{L}_{\ell - 2}(K_{n, m})| = nm[(n - 1)(m - 1)]^{\frac{\ell - 3}{2}}$ for some $\ell \geqslant 3$. For each $R \in \mathscr{L}_{\ell - 2}(K_{n, m})$, we have $|\mathscr{L}_{\ell}(R)| = (n - 1)(m - 1)$. So $|\mathscr{L}_\ell(K_{n, m})| = |\mathscr{L}_{\ell - 2}(K_{n, m})|(n - 1)(m - 1) = nm[(n - 1)(m - 1)]^{\frac{\ell - 1}{2}}$ as desired. The even $\ell$ case is similar.
\qed

\subsection{Loops and multiplicity}
Our next observation is a prerequisite for the study of the chromatic number since it indicates that $\ell$-link graphs are loopless.
\begin{obs}\label{obs_loopless}
For each $(\ell + 1)$-arc $\vec{Q}$, we have $\vec{Q}[0, \ell] \neq \vec{Q}[1, \ell + 1]$.
\end{obs}
\pf Let $G$ be a graph, and $\vec{Q} := (v_0, e_1, \ldots, e_{\ell + 1}, v_{\ell + 1}) \in \vec{\mathscr{L}}_{\ell + 1}(G)$. Since $G$ is loopless, $v_0 \neq v_1$ and hence $\vec{Q}(0, \ell) \neq \vec{Q}(1, \ell + 1)$. So the statement holds for $\ell = 0$. Now let $\ell \geqslant 1$. Suppose for a contradiction that $\vec{Q}(0, \ell) = -\vec{Q}(1, \ell + 1)$. Then $v_i = v_{\ell + 1 - i}$ and $e_{i + 1} = e_{\ell + 1 - i}$ for $i \in \{0, 1, \ldots, \ell\}$. If $\ell = 2s$ for some integer $s \geqslant 1$, then $v_s = v_{s + 1}$, contradicting that $G$ is loopless.  If $\ell = 2s + 1$ for some $s \geqslant 0$, then $e_{s + 1} = e_{s + 2}$, contradicting the definition of a $*$-arc.
\qed



The following statement indicates that, for each $\ell \geqslant 1$, $\mathbb{L}_\ell(G)$ is simple if $G$ is simple, and has multiplicity exactly $2$ otherwise.
\begin{obs}\label{obs_mul2}
Let $G$ be a graph, $\ell \geqslant 1$, and $L_0, L_1 \in \mathscr{L}_{\ell}(G)$. Then $L_0$ can be shunted to $L_1$ through two $(\ell + 1)$-links of $G$ if and only if $G$ contains a $2$-cycle $O := [v_0, e_0, v_1, e_1, v_0]$, such that one of the following cases holds:
\begin{itemize}
\item[{\bf (1)}] $\ell \geqslant 1$ is odd, and $L_i = [v_i, e_i, v_{1 - i}, e_{1 - i}, \ldots, v_i, e_i, v_{1 - i}] \in \mathscr{L}_\ell(O)$ for $i \in \{0, 1\}$. In this case, $[v_i, e_i, v_{1 - i}, e_{1 - i}, \ldots, v_{1 - i}, e_{1 - i}, v_i] \in \mathscr{L}_{\ell + 1}(O)$, for $i \in \{0, 1\}$, are the only two $(\ell + 1)$-links available for the shunting.

\item[{\bf (2)}] $\ell \geqslant 2$ is even, and $L_i = [v_i, e_i, v_{1 - i}, e_{1 - i}, \ldots, v_{1 - i}, e_{1 - i}, v_i] \in \mathscr{L}_\ell(O)$ for $i \in \{0, 1\}$. In this case, $[v_i, e_i, v_{1 - i}, e_{1 - i}, \ldots, v_i, e_i, v_{1 - i}] \in \mathscr{L}_{\ell + 1}(O)$, for $i \in \{0, 1\}$, are the only two $(\ell + 1)$-links available for the shunting.
\end{itemize}
\end{obs}
\pf
$(\Leftarrow)$ is trivial. For $(\Rightarrow)$, since $L_0$ can be shunted to $L_1$, there exists $\vec{L} := (v_0, e_0, v_1, \ldots, v_\ell, e_\ell, v_{\ell + 1}) \in \vec{\mathscr{L}}_{\ell + 1}(G)$ such that $L_i = \vec{L}[i, \ell + i]$ for $i \in \{0, 1\}$. Let $\vec{R} \in \vec{\mathscr{L}}_{\ell + 1}(G) \setminus \{\vec{L}\}$ such that $L_i = \vec{R}[i, \ell + i]$. Then $\vec{L}(i, \ell + i)$ equals $\vec{R}(i, \ell + i)$ or $\vec{R}(\ell + i, i)$. Suppose for a contradiction that $\vec{L}(0, \ell) = \vec{R}(0, \ell)$. Then $\vec{L}(1, \ell) = \vec{R}(1, \ell)$. Since $\vec{L} \neq \vec{R}$, we have $\vec{L}(1, \ell + 1) \neq \vec{R}(1, \ell + 1)$. Thus $\vec{L}(1, \ell + 1) = \vec{R}(\ell + 1, 1)$, and hence $\vec{L}(2, \ell + 1) = \vec{R}(\ell, 1) = \vec{L}(\ell, 1)$, contradicting Observation \ref{obs_loopless}. So $\vec{L}(0, \ell) = \vec{R}(\ell, 0)$. Similarly, $\vec{L}(1, \ell) = \vec{R}(\ell + 1, 1)$. Consequently, $\vec{L}(0, \ell - 1) = \vec{R}(\ell, 1) = \vec{L}(2, \ell + 1)$; that is, $v_j = v_0$ and $e_j = e_0$ if $j \in [0, \ell]$ is even, while $v_j = v_1$ and $e_j = e_1$ if $j \in [0, \ell + 1]$ is odd.
\qed


\subsection{Connectedness} This subsection characterises when $\mathbb{L}_\ell(G)$ is connected.
A middle segment of $L \in \mathscr{L}_\ell(G)$ is a \emph{middle unit}, written \emph{$c_L$}, if it is a unit of $G$. Note that $c_L$ is a vertex if $\ell$ is even, and is an edge otherwise. Denote by \emph{$G(\ell)$} the subgraph of $G$ induced by the middle units of $\ell$-links of $G$.



The lemma below is important in dealing with the connectedness of $\ell$-link graphs. Before stating it, we define a \emph{conjunction} operation, which is an extension of an operation by Biggs \cite[Chapter 17]{Biggs1993}. Let $\vec{L} := (v_0, e_1, v_1, \ldots, e_\ell, v_\ell) \in \vec{\mathscr{L}}_\ell(G)$ and $\vec{R} := (u_0, f_1, u_1, \ldots, f_s, u_s) \in \vec{\mathscr{L}}_s(G)$ such that $v_\ell = u_0$ and $e_{\ell} \neq f_1$. The \emph{conjunction} of  $\vec{L}$ and $\vec{R}$ is $(\vec{L}.\vec{R}) := (v_0, e_1, \ldots, e_{\ell}, v_{\ell} = u_0, f_1, \ldots, f_s, u_s) \in \vec{\mathscr{L}}_{\ell + s}(G)$ or $[\vec{L}.\vec{R}] := [v_0, e_1, \ldots, e_{\ell}, v_{\ell} = u_0, f_1, \ldots, f_s, u_s] \in \mathscr{L}_{\ell + s}(G)$.

\begin{lem}\label{lem_GlConnect}
Let $\ell, s \geqslant 0$, and $G$ be a connected graph. Then $G(\ell)$ is connected. And each $s$-link of $G(\ell)$ is a middle segment of a $(2 \lfloor\frac{\ell}{2}\rfloor + s)$-link of $G$. Moreover, for $\ell$-links $L$ and $R$ of $G$, there is an $\ell$-link $L'$ with middle unit $c_L$, and an $\ell$-link $R'$ with middle unit $c_R$, such that $L'$ can be shunted to $R'$.
\end{lem}
\pf For $\ell \in \{0, 1\}$, since $G$ is connected, $G(\ell) = G$ and the lemma holds. Let $\ell := 2m \geqslant 2$ be even. $u, v \in V(G(\ell))$ if and only if they are middle vertices of some $\vec{L}, \vec{R} \in \vec{\mathscr {L}}_\ell(G)$ respectively. Since $G$ is connected, there exists some $\vec{P} \in \vec{\mathscr{L}}_s(G)$ from $(u, e, u_1)$ to $(v_{s - 1}, f, v)$. By Observation \ref{obs_loopless}, $\vec{L}[m - 1, m] \neq \vec{L}[m, m + 1]$. For such an $s$-arc $\vec{P}$, without loss of generality, $e \neq \vec{L}[m - 1, m]$, and similarly, $f \neq \vec{R}[m, m + 1]$. Then $\vec{P}$ is a middle segment of $\vec{Q} := (\vec{L}(0, m).\vec{P}.\vec{R}(m, 2m)) \in \vec{\mathscr{L}}_{\ell + s}(G)$. So $\vec{P} \in \vec{\mathscr{L}}_{s}(G(\ell))$. And $L' := \vec{Q}[0, \ell]$ can be shunted to $R' := \vec{Q}[s, \ell + s]$ through $\vec{Q}$. The odd $\ell$ case is similar.
\qed




Sufficient conditions for $\mathbb{A}_\ell(G)$ to be strongly connected
can be found in \cite[Page 76]{GR2001}. The following corollary of Lemma \ref{lem_GlConnect} reveals a strong relationship between the shunting of $\ell$-links and the connectedness of $\ell$-link graphs.
\begin{coro}\label{coro_LlConnect}
For a connected graph $G$, $\mathbb{L}_\ell(G)$ is connected if and only if any two $\ell$-links of $G$ with the same middle unit can be shunted to each other.
\end{coro}

We now present our main result of this section, which plays a key role in dealing with the graph minors of $\ell$-link graphs in Section \ref{sec_Hadwiger}.
\begin{lem}\label{lem_hub}
Let $G$ be a graph, and $X$ be a connected subgraph of $G(\ell)$. Then for every pair of $\ell$-links $L$ and $R$ of $X$, $L$ can be shunted to $R$ under the restriction that in each step, the middle unit of the image of $L$ belongs to $X$.
\end{lem}
\pf
First we consider the case that $c_L$ is in $R$. Then there is a common segment $Q$ of $L$ and $R$ of maximum length containing $c_L$. Without loss of generality, assign directions to $L$ and $R$ such that $\vec{L} = (\vec{L}_0.\vec{Q}.\vec{L}_1)$ and $\vec{R} = (\vec{R}_1.\vec{Q}. \vec{R}_0)$, where $\vec{L}_i \in \vec{\mathscr{L}}_{\ell_i}(X)$ and $\vec{R}_i \in \vec{\mathscr{L}}_{s_i}(X)$ for $i \in \{0, 1\}$ such that $s_1 \geqslant s_0$. Then $\ell \geqslant \ell_0 + \ell_1 = s_0 + s_1 \geqslant s_1$. Let $x$ be the head vertex and $e$  be the head edge of $\vec{L}$. Since $c_L$ is in $Q$, $\ell_0 \leqslant \ell/2$. Since $X$ is a subgraph of $G(\ell)$, by Lemma \ref{lem_GlConnect}, there exists $\vec{L}_2 \in \vec{\mathscr{L}}_{\ell_0}(G)$ with tail vertex $x$ and tail edge different from $e$. Let $y$ be the tail vertex and $f$ be the tail edge of $\vec{R}$. Then there exits $\vec{R}_2 \in \vec{\mathscr{L}}_{s_0}(G)$ with head vertex $y$ and head edge different from $f$. We can shunt $L$ to $R$ first through $(\vec{L}.\vec{L}_2) \in \vec{\mathscr{L}}_{\ell + \ell_0}(G)$, then $-(\vec{R}_2.\vec{R}_1.\vec{Q}.\vec{L}_1.\vec{L}_2) \in \vec{\mathscr{L}}_{\ell + \ell_0 + \ell_1}(G)$, and finally $(\vec{R}_2.\vec{R}) \in \vec{\mathscr{L}}_{\ell + s_0}(G)$. Since $\ell_0 \leqslant \ell/2$ and $s_0 \leqslant s_1 \leqslant \ell/2$, the middle unit of each image is inside $L$ or $R$.

Secondly, we consider the case that $c_L$ is not in $R$. Then there exists a segment $Q$ of $L$ of maximum length that contains $c_L$, and is edge-disjoint with $R$. Since $X$ is connected, there exists a shortest $*$-arc $\vec{P}$ from a vertex $v$ of $R$ to a vertex $u$ of $L$. Then $P$ is edge-disjoint with $Q$ because of its minimality. Without loss of generality, assign directions to $L$ and $R$ such that $u$ separates $\vec{L}$ into $(\vec{L}_0.\vec{L}_1)$ with $c_L$ on $L_1$, and $v$ separates $\vec{R}$ into $(\vec{R}_1.\vec{R}_0)$, where $L_i$ is of length $\ell_i$ while $R_i$ is of length $s_i$ for $i \in \{0, 1\}$, such that $s_1 \geqslant s_0$. Then $\ell_0, s_0 \leqslant \ell/2$. Let $x$ be the head vertex and $e$ be the head edge of $\vec{L}$. Since $\ell_0 \leqslant \ell/2$ and $X$ is a subgraph of $G(\ell)$, by Lemma \ref{lem_GlConnect}, there exists an $\ell_0$-arc $\vec{L}_2$ of $G$ with tail vertex $x$ and tail edge different from $e$. Let $y$ be the tail vertex and $f$ be the tail edge of $\vec{R}$. Then there exits an $s_0$-arc $\vec{R}_2$ of $G$ with head vertex $y$ and head edge different from $f$. Now we can shunt $L$ to $R$ through $(\vec{L}.\vec{L}_2)$, $-(\vec{R}_2.\vec{R}_1.\vec{P}.\vec{L}_1.\vec{L}_2)$ and $(\vec{R}_2.\vec{R})$ consecutively. One can check that in this process the middle unit of each image belongs to $L, P$ or $R$.
\qed

From Lemma \ref{lem_hub}, the set of $\ell$-links of a connected $G(\ell)$ serves as a `hub' in the shunting of $\ell$-links of $G$. More explicitly, for $L, R \in \mathscr{L}_\ell(G)$, if we can shunt $L$ to $L' \in \mathscr{L}_\ell(G(\ell))$, and $R$ to $R' \in \mathscr{L}_\ell(G(\ell))$, then $L$ can be shunted to $R$ since $L'$ can be shunted to $R'$. Thus we have the following corollary which provides a more efficient way to test the connectedness of $\ell$-link graphs.
\begin{coro}\label{coro_hub} Let $G$ be a graph. Then
$\mathbb{L}_\ell(G)$ is connected if and only if $G(\ell)$ is connected, and each $\ell$-link of $G$ can be shunted to an $\ell$-link of $G(\ell)$.
\end{coro}

%% file: Sec35.tex
\section{Chromatic number of $\ell$-link graphs}\label{sec_chromatic}
In this section, we reveal a recursive structure of $\ell$-link graphs, which leads to an upper bound for the chromatic number of $\ell$-link graphs.

\begin{lem}\label{lem_VEpartition}
Let $G$ be a graph and $\ell \geqslant 2$ be an integer. Then
$(\mathcal {V}, \mathcal {E}) := (\mathcal {V}_{\ell}(G), \mathcal {E}_{\ell}(G))$ is an almost standard partition of $H := \mathbb{L}_\ell(G)$. Further, $H_{(\mathcal {V}, \mathcal {E})}$ is isomorphic to an induced subgraph of $\mathbb{L}_{\ell - 2}(G)$.
\end{lem}
\pf First we verify that $(\mathcal {V}, \mathcal {E})$ is an almost standard partition of $H$.

{\bf (a)} We prove that, for each $R \in \mathscr{L}_{\ell - 2}(G)$, $V := \mathscr{L}_\ell(R) \in \mathcal{V}$ is an independent set of $H$. Suppose not. Then there are $\vec{L}, \vec{L}' \in \vec{\mathscr{L}}_\ell(G)$ such that $L, L' \in V$, and $L$ can be shunted to $L'$ in one step. Then $R = \vec{L}[1, \ell - 1]$ can be shunted to $R = \vec{L}'[1, \ell - 1]$ in one step, contradicting Observation \ref{obs_loopless}.

{\bf (b)} Here we show that each $E \in \mathcal {E}$ is incident to exactly two parts of $\mathcal {V}$. By definition there exists $P \in \mathscr {L}_{\ell - 1}(G)$ with $\mathscr{L}_{\ell + 1}(P) = E$. Let $\{L, R\} := P^{\{\ell - 2\}}$. Then $\mathscr{L}_{\ell}(L)$ and $\mathscr{L}_{\ell}(R)$ are the only two parts of $\mathcal {V}$ incident to $E$.

{\bf (c)} We explain that each $E \in \mathcal{E}$ is the edge set of a complete bipartite subgraph of $H$. By definition there exists $\vec{P} \in \vec{\mathscr {L}}_{\ell - 1}(G)$ with $\mathscr{L}_{\ell + 1}(P) = E$. Let $A := \{[\vec{e}.\vec{P}] \in \mathscr{L}_\ell(G)\}$ and $B := \{[\vec{P}.\vec{f}] \in \mathscr{L}_\ell(G)\}$. One can check that $E$ induces a complete bipartite subgraph of $H$ with bipartition $A \cup B$.

{\bf (d)} We prove that each $v \in V(H)$ is incident to at most two parts of $\mathcal{E}$. By definition there exists $Q \in \mathscr{L}_\ell(G)$ with $Q = v$. Then the set of edge parts of $\mathcal{E}$ incident to $v$ is $\{\mathscr {L}_{\ell + 1}(L) \neq \emptyset| L \in Q^{\{\ell - 1\}}\}$ with cardinality at most $2$.

{\bf (e)} Let $v$ be a vertex of $V \in \mathcal{V}$ incident to different $E, F \in \mathcal {E}$. We explain that $v$ is uniquely determined by $V$, $E$ and $F$.
By definition there exists $\vec{P} \in \vec{\mathscr {L}}_{\ell - 2}(G)$ such that $V = \mathscr {L}_{\ell}(P)$. There also exists $Q := [\vec{e}_1.\vec{P}.\vec{e}_{\ell}] \in \mathscr {L}_{\ell}(P)$ such that $v = Q$. Besides, there are $L, R \in \mathscr{L}_{\ell - 1}(G)$ such that $E = \mathscr{L}_{\ell + 1}(L)$ and $F = \mathscr{L}_{\ell + 1}(R)$. Then $\{L, R\} = Q^{\{\ell - 1\}}$ since $L \neq R$. Note that $Q$ is uniquely determined by $Q^{\{\ell - 1\}}$ and $c_Q = c_P$. Thus it is uniquely determined by $E = \mathscr{L}_{\ell + 1}(L), F = \mathscr{L}_{\ell + 1}(R)$ and $V = \mathscr {L}_{\ell}(P)$.

Now we show that $H_{(\mathcal {V}, \mathcal {E})}$ is isomorphic to an induced subgraph of $\mathbb{L}_{\ell - 2}(G)$. Let $X$ be the subgraph of $\mathbb{L}_{\ell - 2}(G)$ of vertices $L \in \mathscr{L}_{\ell - 2}(G)$ such that
$\mathscr {L}_{\ell}(L) \neq \emptyset$, and edges $Q \in \mathscr {L}_{\ell - 1}(G)$ such that $\mathscr {L}_{\ell + 1}(Q) \neq \emptyset$. One can check that $X$ is an induced subgraph of $\mathbb{L}_{\ell - 2}(G)$. An isomorphism from $H_{(\mathcal {V}, \mathcal {E})}$ to $X$ can be defined as the injection sending $\mathscr {L}_{\ell}(L) \neq \emptyset$ to $L$, and $\mathscr {L}_{\ell + 1}(Q) \neq \emptyset$ to $Q$.
\qed

Below we give an interesting algorithm for colouring a class of graphs.
\begin{lem}\label{lem_colourr}
Let $H$ be a graph with a $t$-colouring such that each vertex of $H$ is adjacent to at most $r \geqslant 0$ differently coloured vertices. Then $\chi(H) \leqslant \lfloor \frac{tr}{r + 1} \rfloor + 1$.
\end{lem}
\pf
The result is trivial for $t = 0$ since, in this case, $\chi(H) = 0$. If $r + 1 \geqslant t \geqslant 1$, then $\lfloor \frac{tr}{r + 1} \rfloor + 1 = t$, and the lemma holds since $t \geqslant \chi(H)$.


Now assume $t \geqslant r + 2 \geqslant 2$. Let $U_1, U_2, \ldots, U_t$ be the colour classes of the given colouring. For $i \in [t]$, denote by $i$ the colour assigned to vertices in $U_i$. Run the following algorithm: For $j = 1, \ldots, t$, and for each $u \in U_{t - j + 1}$, let $s \in [t]$ be the minimum integer that is not the colour of a neighbour of $u$ in $H$; if $s < t - j + 1$, then recolour $u$ by $s$.

In the algorithm above, denote by $C_i$ the set of colours used by the vertices in $U_i$ for $i \in [t]$. Let $k := \lfloor\frac{t - 1}{r + 1}\rfloor$. Then $t - 1 \geqslant k(r + 1) \geqslant k \geqslant 1$. We claim that after $j \in [0, k]$ steps, $C_{t - i + 1} \subseteq [ir + 1]$ for $i \in [j]$, and $C_i = \{i\}$ for $i \in [t - j]$. This is trivial for $j = 0$. Inductively assume it holds for some $j \in [0, k - 1]$. In the $(j + 1)^{th}$ step, we change the colour of each $u \in U_{t - j}$ from $t - j$ to the minimum $s \in [t]$ that is not used by the neighbourhood of $u$. It is enough to show that $s \leqslant (j + 1)r + 1$.

First suppose that all neighbours of $u$ are in $\bigcup_{i \in [t - j - 1]}U_i$. By the analysis above, $t - j - 1 \geqslant t - k \geqslant kr + 1 \geqslant r + 1$. So at least one part of $\mathcal{S} := \{U_i| i \in [t - j - 1]\}$ contains no neighbour of $u$. From the induction hypothesis, $C_i = \{i\}$ for $i \in [t - j - 1]$. Hence at least one colour in $[r + 1]$ is not used by the neighbourhood of $u$; that is, $s \leqslant r + 1 \leqslant (j + 1)r + 1$.

Now suppose that $u$ has at least one neighbour in $\bigcup_{i \in [t - j + 1, t]}U_i$. By the induction hypothesis, $\bigcup_{i \in [t - j + 1, t]}C_i \subseteq [jr + 1]$. At the same time, $u$ has neighbours in at most $r - 1$ parts of $\mathcal{S}$. So the colours possessed by the neighbourhood of $u$ are contained in $[jr + 1 + r - 1] = [(j + 1)r]$. Thus $s \leqslant (j + 1)r + 1$. This proves our claim.

The claim above indicates that, after the $k^{th}$ step, $C_{t - i + 1} \subseteq [ir + 1]$ for $i \in [k]$, and $C_i = \{i\}$ for $i \in [t - k]$. Hence we have a $(t - k)$-colouring of $H$ since $t - k \geqslant kr + 1$. Therefore, $\chi(H) \leqslant t - k = \lceil\frac{tr + 1}{r + 1}\rceil = \lfloor \frac{tr}{r + 1} \rfloor + 1$.
\qed

Lemma \ref{lem_VEpartition} indicates that $\mathbb{L}_\ell(G)$ is homomorphic to $\mathbb{L}_{\ell - 2}(G)$ for $\ell \geqslant 2$. So by \cite[Proposition 1.1]{Cameron2006}, $\chi_\ell(G) \leqslant \chi_{\ell - 2}(G)$. By Lemma \ref{lem_VEpartition}, every vertex of $\mathbb{L}_\ell(G)$ has neighbours in at most two parts of $\mathcal {V}_\ell(G)$, which enables us to improve the upper bound on $\chi_\ell(G)$.

\begin{lem}\label{lem_colourr}
Let $G$ be a graph, and $\ell \geqslant 2$. Then $\chi_\ell(G) \leqslant \lfloor \frac{2}{3}\chi_{\ell - 2}(G)\rfloor + 1$.
\end{lem}
\pf
By Lemma \ref{lem_VEpartition}, $(\mathcal {V}, \mathcal {E}) := (\mathcal {V}_{\ell}(G), \mathcal {E}_{\ell}(G))$ is an almost standard partition of $H := \mathbb{L}_\ell(G)$. So each vertex of $H$ has neighbours in at most two parts of $\mathcal{V}$. Further, $H_{\mathcal{V}}$ is a subgraph of $\mathbb{L}_{\ell - 2}(G)$. So $\chi_\ell(G) \leqslant \chi := \chi(H_{\mathcal{V}}) \leqslant \chi_{\ell - 2}(G)$.

We now construct a $\chi$-colouring of $H$ such that each vertex of $H$ is adjacent to at most two differently coloured vertices. By definition $H_{\mathcal{V}}$ admits a $\chi$-colouring with colour classes $K_1, \ldots, K_{\chi}$. For $i \in [\chi]$, assign the colour $i$ to each vertex of $H$ in $U_i := \bigcup_{V_{\mathcal {V}} \in K_i}V$. One can check that this is a desired colouring. In Lemma \ref{lem_colourr}, letting $t = \chi$ and $r = 2$ yields that $\chi_\ell(G) \leqslant \lfloor\frac{2}{3}\chi\rfloor + 1$. Recall that $\chi \leqslant \chi_{\ell - 2}(G)$. Thus the lemma follows. \qed

%
%

As shown below, Lemma \ref{lem_colourr} can be applied recursively to produce an upper bound for $\chi_\ell(G)$ in terms of $\chi(G)$ or $\chi'(G)$.

\medskip
\noindent{\textbf{Proof of Theorem \ref{thm_chrl}}}. When $\ell \in \{0, 1\}$, it is trivial for {\bf (1)}{\bf (2)} and {\bf (4)}.
By \cite[Proposition 5.2.2]{Diestel2010}, $\chi_0 = \chi \leqslant \Delta + 1$. So {\bf (3)} holds. Now let $\ell \geqslant 2$. By Lemma \ref{lem_VEpartition}, $H := \mathbb{L}_\ell(G)$ admits an almost standard partition $(\mathcal{V}, \mathcal{E}) := (\mathcal{V}_\ell(G), \mathcal{E}_\ell(G))$, such that $H_{(\mathcal{V}, \mathcal{E})}$ is an induced subgraph of $\mathbb{L}_{\ell - 2}(G)$. By definition each part of $\mathcal{V}$ is an independent set of $H$. So $H \rightarrow \mathbb{L}_{\ell - 2}(G)$, and $\chi_\ell \leqslant \chi_{\ell - 2}$. This proves {\bf (4)}. Moreover, each vertex of $H$ has neighbours in at most two parts of $\mathcal{V}$. By Lemma \ref{lem_colourr}, $\chi_{\ell} := \chi_\ell(G) \leqslant \frac{2\chi_{\ell - 2}}{3} + 1$. Continue the analysis, we have
$\chi_\ell \leqslant \chi_{\ell - 2i}$, and $\chi_\ell - 3 \leqslant (\frac{2}{3})^{i}(\chi_{\ell - 2i} - 3)$ for $1 \leqslant i \leqslant \lfloor \ell/2\rfloor$. Therefore, if $\ell$ is even, then $\chi_\ell \leqslant \chi_0 = \chi \leqslant \Delta + 1$, and $\chi_\ell - 3 \leqslant  (\frac{2}{3})^{\ell/2}(\chi - 3)$. Thus {\bf (1)} holds. Now let $\ell \geqslant 3$ be odd. Then $\chi_\ell \leqslant \chi_1 = \chi'$, and $\chi_\ell - 3 \leqslant (\frac{2}{3})^{\frac{\ell - 1}{2}}(\chi' - 3)$. This verifies {\bf (2)}. As a consequence, $\chi_\ell \leqslant \chi_3 \leqslant \frac{2}{3}(\chi' - 3) + 3 = \frac{2}{3}\chi' + 1$. By Shannon \cite{Shannon1949}, $\chi' \leqslant \frac{3}{2}\Delta$. So $\chi_\ell \leqslant \Delta + 1$, and hence {\bf (3)} holds.
\qed

%
%

The following corollary of Theorem \ref{thm_chrl} implies that Hadwiger's conjecture is true for $\mathbb{L}_\ell(G)$ if $G$ is regular and $\ell \geqslant 4$.
\begin{coro}\label{coro_chrHadl}
Let $G$ be a graph with $\Delta := \Delta(G) \geqslant 3$. Then $\chi_\ell(G) \leqslant 3$ for all $\ell > 2\log_{1.5}(\Delta - 2) + 3$. Further, Hadwiger's conjecture holds for $\mathbb{L}_\ell(G)$ if $\ell > 2\log_{1.5}(\Delta - 2) - 3.83$, or $\dege := \dege(G) \geqslant 3$ and $\ell > 2\log_{1.5}\frac{\Delta - 2}{\dege - 2} + 3$.
\end{coro}
\pf
By Theorem \ref{thm_chrl}, for each $t \geqslant 3$, $\chi_\ell := \chi_\ell(G) \leqslant t$ if $(\frac{2}{3})^{\ell/2}(\Delta - 2)  < t - 2$ and $(\frac{2}{3})^{\frac{\ell - 1}{2}}(\frac{3}{2}\Delta - 3) < t - 2$. Solving these inequalities gives $\ell > 2\log_{1.5}(\Delta - 2) - 2\log_{1.5}(t - 2) + 3$. Thus $\chi_\ell \leqslant 3$ if $\ell > 2\log_{1.5}(\Delta - 2) + 3$. So the first statement holds. By Robertson et al. \cite{RST1993} and Theorem \ref{thm_hgeqg}, Hadwiger's conjecture holds for $\mathbb{L}_\ell(G)$ if $\ell \geqslant 1$ and $\chi_\ell \leqslant \max\{6, \dege\}$. Letting $t = 6$ gives that $\ell > 2\log_{1.5}(\Delta - 2) - 4\log_{1.5}2 + 3$. Letting $t = \dege \geqslant 3$ gives that $\ell > 2\log_{1.5}\frac{\Delta - 2}{\dege - 2} + 3$. So the corollary holds since $4\log_{1.5}2 - 3 > 3.83$.
\qed

\noindent{\textbf{Proof of Theorem \ref{thm_hcllink}(3)(4)(5)}}. {\bf (3)} and {\bf (4)} follow from Corollary \ref{coro_chrHadl}. Now consider {\bf (5)}. By Reed and Seymour \cite{Reed2004}, Hadwiger's conjecture holds for $\mathbb{L}_1(G)$. If $\ell \geqslant 2$ and $\Delta \leqslant 5$, by Theorem \ref{thm_chrl}(3), $\chi_\ell(G) \leqslant 6$. In this case, Hadwiger's conjecture holds for $\mathbb{L}_\ell(G)$ by Robertson et al.  \cite{RST1993}.\qed

%% file: Sec36.tex
\section{Complete minors of $\ell$-link graphs}\label{sec_Hadwiger}
It has been proved in the last section that Hadwiger's conjecture is true for $\mathbb{L}_{\ell}(G)$ if $\ell$ is large enough. In this section, we further investigate the minors, especially the complete minors, of $\ell$-link graphs. To see the intuition of our method, let $v$ be a vertex of degree $t$ in $G$. Then $\mathbb{L}_1(G)$ contains a $K_t$-subgraph whose vertices correspond to the edges of $G$ incident to $v$. For $\ell \geqslant 2$, roughly speaking, we extend $v$ to a subgraph $X$ of diameter less than $\ell$, and extend each edge incident to $v$ to an $\ell$-link of $G$ starting from a vertex of $X$. By studying the shunting of these $\ell$-links, we find a $K_t$-minor in $\mathbb{L}_\ell(G)$.

For subgraphs $X, Y$ of $G$, let \emph{$\vec{E}(X, Y)$} be the set of arcs of $G$ from $V(X)$ to $V(Y)$, and \emph{$E(X, Y)$} be the set of edges of $G$ between $V(X)$ and $V(Y)$.

\begin{lem}\label{lem_ltConnected}
Let $\ell \geqslant 1$ be an integer, $G$ be a graph, and $X$ be a subgraph of $G$ with $\diam(X) < \ell$ such that $Y := G - V(X)$ is connected. If $t := |E(X, Y)| \geqslant 2$, then $\mathbb{L}_\ell(G)$ contains a $K_t$-minor.
\end{lem}
\pf Let $\vec{e}_1, \ldots, \vec{e}_t$ be distinct arcs in $\vec{E}(Y, X)$. Say $\vec{e}_i = (y_i, e_i, x_i)$ for $i \in [t]$. Since $\diam(X) < \ell$, there is a dipath  $\vec{P}_{ij}$ of $X$ from $x_i$ to $x_j$ of length $\ell_{ij} \leqslant \ell - 1$ such that $P_{ij} = P_{ji}$. Since $Y$ is connected, it contains a  dipath $\vec{Q}_{ij}$ from $y_i$ to $y_j$. Since $t \geqslant 2$, $O_i := [\vec{P}_{i\; i'}.-\vec{e}_{i'}.\vec{Q}_{i'\; i}. \vec{e}_i]$ is a cycle of $G$, where $i' := (i \mod t) + 1$. Thus $H := \mathbb{L}_\ell(G)$ contains a cycle $\mathbb{L}_\ell(O_1)$, and hence a $K_2$-minor. Now let $t \geqslant 3$, and $\vec{L}_i \in \vec{\mathscr{L}}_\ell(O_i)$ with head arc $\vec{e}_i$. Then $[\vec{L}_i.\vec{P}_{ij}]^{[\ell]} \in \mathscr{L}_{\ell_{ij}}(H)$. And the union of the units of $[\vec{L}_i.\vec{P}_{ij}]^{[\ell]}$ over $j \in [t]$ is a connected subgraph $X_i$ of $H$. In the remainder of the proof, for distinct $i, j \in [t]$, we show that $X_i$ and $X_j$ are disjoint. Further, we construct a path in $H$ between $X_i$ and $X_j$ that is internally disjoint with its counterparts, and has no inner vertex in any of $V(X_1), \ldots, V(X_t)$. Then by contracting each $X_i$ into a vertex, and each path into an edge, we obtain a $K_t$-minor of $H$.

First of all, assume for a contradiction that there are different $i, j \in [t]$ such that $X_i$ and $X_j$ share a common vertex that corresponds to an $\ell$-link $R$ of $G$. Then by definition, there exists some $p \in [t]$ such that $R$ can be obtained by shunting $L_i$ along $(\vec{L}_i.\vec{P}_{ip})$ by some $s_i \leqslant \ell_{ip}$ steps. So $R = [\vec{L}_i(s_i, \ell).\vec{P}_{ip}(0, s_i)]$. Similarly, there are $q \in [t]$ and $s_j \leqslant \ell_{jq}$ such that $R = [\vec{L}_j(s_j, \ell).\vec{P}_{jq}(0, s_j)]$. Recall that $E(X) \cap E(X, Y) = E(Y) \cap E(X, Y) = \emptyset$. So $e_i = \vec{L}_i[\ell - 1, \ell]$ and $e_j = \vec{L}_j[\ell - 1, \ell]$ belong to both $L_i$ and $L_j$. By the definition of $O_i$, this happens if and only if $i = j'$ and $j = i'$, which is impossible since $t \geqslant 3$.

Secondly, for different $i, j \in [t]$, we define a path of $H$ between $X_i$ and $X_j$. Clearly, $L_i$ can be shunted to $L_j$ through $\vec{R}_{ij}' := (\vec{L}_i.\vec{P}_{ij}.-\vec{L}_j)$ in $G$. In this shunting, $L_i' := [\vec{L}_i(\ell_{ij}, \ell).\vec{P}_{ij}]$ is the last image corresponding to a vertex of $X_i$, while $L_j' := [\vec{P}_{ij}.\vec{L}_j(\ell, \ell_{ij})]$ is the first image corresponding to a vertex of $X_j$. Further, $L_i'$ can be shunted to $L_j'$ through $\vec{R}_{ij} := (\vec{L}_i(\ell_{ij}, \ell).\vec{P}_{ij}.\vec{L}_j(\ell, \ell_{ij})) \in \vec{\mathscr{L}}_{2\ell - \ell_{ij}}(G)$, which is a subsequence of $\vec{R}_{ij}'$. Then $R_{ij}^{[\ell]}$ is an $(\ell - \ell_{ij})$-path of $H$ between $X_i$ and $X_j$. We show that for each $p \in [t]$, $X_p$ contains no inner vertex of $R_{ij}^{[\ell]}$. When $\ell - \ell_{ij} = 1$, $R_{ij}^{[\ell]}$ contains no inner vertex. Now assume $\ell - \ell_{ij} \geqslant 2$. Each inner vertex of $R_{ij}^{[\ell]}$ corresponds to some $Q_{ij} := [\vec{L}_i(s_i, \ell).\vec{P}_{ij}.\vec{L}_j(\ell, \ell + \ell_{ij} - s_i)] \in \mathscr{L}_\ell(G)$, where $\ell_{ij} + 1 \leqslant s_i \leqslant \ell - 1$. Assume for a contradiction that for some $p \in [t]$, $X_p$ contains a vertex corresponding to $Q_{ij}$. By definition there exists $q \in [t]$ such that $Q_{ij} = [\vec{L}_p(s_p, \ell).\vec{P}_{pq}(0, s_p)]$, where $0 \leqslant s_p \leqslant \ell_{pq}$. Without loss of generality, $(\vec{L}_i(s_i, \ell).\vec{P}_{ij}.\vec{L}_j(\ell, \ell + \ell_{ij} - s_i)) = (\vec{L}_p(s_p, \ell).\vec{P}_{pq}(0, s_p))$. Since $e_j$ and $e_p$ are not in $P_{pq}$, hence $\vec{e}_j$ belongs to $-\vec{L}_p$ and $\vec{e}_p$ belongs to $-\vec{L}_j$. By the definition of $\vec{L}_i$, this happens only when $j = p'$ and $p = j'$, contradicting $t \geqslant 3$.




We now show that $R_{ij}^{[\ell]}$ and $R_{pq}^{[\ell]}$ are internally disjoint, where $i \neq j$, $p \neq q$ and $\{i, j\} \neq \{p, q\}$. Suppose not. Then by the analysis above, there are $s_i$ and $s_p$ with $\ell_{ij} + 1 \leqslant s_i \leqslant \ell - 1$ and $\ell_{pq} + 1 \leqslant s_p \leqslant \ell - 1$ such that $Q_{ij} = Q_{pq}$. Without loss of generality, $(\vec{L}_i(s_i, \ell).\vec{P}_{ij}.\vec{L}_j(\ell, \ell + \ell_{ij} - s_i)) = (\vec{L}_p(s_p, \ell).\vec{P}_{pq}.\vec{L}_q(\ell, \ell + \ell_{pq} - s_p))$. If $s_i = s_p$, then $\vec{e}_i = \vec{e}_p$ and $\vec{e}_j = \vec{e}_q$ since $E(X) \cap E(X, Y) = \emptyset$; that is, $i = p$ and $j = q$, contradicting $\{i, j\} \neq \{p, q\}$. Otherwise, with no loss of generality, $s_i > s_p$. Then $\vec{e}_q$ and $\vec{e}_i$ belong to $\vec{L}_j$ and $\vec{L}_p$ respectively; that is, $i = p$ and $j = q$, again contradicting $\{i, j\} \neq \{p, q\}$.


In summary, $X_1, \ldots, X_t$ are vertex-disjoint connected subgraphs, which are pairwise connected by internally disjoint  $*$-links $R_{ij}^{[\ell]}$ of $H$, such that no inner vertex of $R_{ij}^{[\ell]}$ is in $V(X_1)\cup \ldots \cup V(X_t)$. So by contracting each $X_i$ to a vertex, and $R_{ij}^{[\ell]}$ to an edge, we obtain a $K_t$-minor of $H$. \qed

\begin{lem}\label{lem_ltConCycle}
Let $\ell \geqslant 1$, $G$ be a graph, and $X$ be a subgraph of $G$ with $\diam(X) < \ell$ such that $Y := G - V(X)$ is connected and contains a cycle. Let $t := |E(X, Y)|$. Then $\mathbb{L}_\ell(G)$ contains a $K_{t + 1}$-minor.
\end{lem}
\pf Let $O$ be a cycle of $Y$. Then $H := \mathbb{L}_\ell(G)$ contains a cycle $\mathbb{L}_\ell(O)$ and hence a $K_2$-minor. Now assume $t \geqslant 2$. Let $\vec{e}_1, \ldots, \vec{e}_t$ be distinct arcs in $\vec{E}(Y, X)$. Say $\vec{e}_i = (y_i, e_i, x_i)$ for $i \in [t]$. Since $Y$ is connected, there is a dipath $\vec{P}_i$ of $Y$ of minimum length $s_i \geqslant 0$ from some vertex $z_i$ of $O$ to $y_i$. Let $\vec{Q}_i$ be an $\ell$-arc of $O$ with head vertex $z_i$. Then $\vec{L}_i := (\vec{Q}_i.\vec{P}_i.\vec{e}_i)(s_i + 1, \ell + s_i + 1) \in \vec{\mathscr{L}}_\ell(G)$. Since $\diam(X) \leqslant \ell - 1$, there is a dipath $\vec{P}_{ij}$ of $X$ of length $\ell_{ij} \leqslant \ell - 1$ from $x_i$ to $x_j$ such that $P_{ij} = P_{ji}$.

Clearly, $[\vec{L}_i.\vec{P}_{ij}]^{[\ell]}$ is an $\ell_{ij}$-link of $H$. And the union of the units of $[\vec{L}_i.\vec{P}_{ij}]^{[\ell]}$ over $j \in [t]$ induces a connected subgraph $X_i$ of $H$. For different $i, j \in [t]$, let $R_{ij} := [\vec{L}_i(\ell_{ij}, \ell).\vec{P}_{ij}.\vec{L}_j(\ell, \ell_{ij})] = R_{ji} \in \mathscr{L}_{2\ell - \ell_{ij}}(G)$. Then $R_{ij}^{[\ell]}$ is an $(\ell - \ell_{ij})$-path of $H$ between $X_i$ and $X_j$. As in the proof of Lemma \ref{lem_ltConnected}, it is easy to check that $X_1, \ldots, X_t$ are vertex-disjoint connected subgraphs of $H$, which are pairwise connected by internally disjoint paths $R_{ij}^{[\ell]}$. Further, no inner vertex of $R_{ij}^{[\ell]}$ is in $V(X_1)\cup \ldots \cup V(X_t)$. So a $K_t$-minor of $H$ is obtained accordingly.

Finally, let $Z$ be the connected subgraph of $H$ induced by the units of $\mathbb{L}_\ell(O)$ and $[\vec{Q}_i.\vec{P}_i]^{[\ell]}$ over $i \in [t]$. Then $Z$ is vertex-disjoint with $X_i$ and with the paths $R_{ij}^{[\ell]}$. Moreover, $Z$ sends an edge $(\vec{Q}_i.\vec{P}_i.\vec{e}_i)(s_i, \ell + s_i + 1)^{[\ell]}$ to each $X_i$. Thus $H$ contains a $K_{t+1}$-minor. \qed

In the following, we use the `hub' (described after Lemma \ref{lem_hub}) to construct certain minors in $\ell$-link graphs.
\begin{coro}\label{coro_lminor}
Let $\ell \geqslant 0$, $G$ be a graph, $M$ be a minor of $G(\ell)$ such that each branch set contains an $\ell$-link. Then $\mathbb{L}_{\ell}(G)$ contains an $M$-minor.
\end{coro}
\pf Let $X_1, \ldots, X_t$ be the branch sets of an $M$-minor of $G(\ell)$ such that $X_i$ contains an $\ell$-link for each $i \in [t]$.
For any connected subgraph $Y$ of $G(\ell)$ contains at least one $\ell$-link, let $\mathbb{L}_{\ell}(G, Y)$ be the subgraph of $H := \mathbb{L}_\ell(G)$ induced by the $\ell$-links of $G$ of which the middle units are in $Y$. Let $H(Y)$ be the union of the components of $\mathbb{L}_\ell(G, Y)$ which contains at least one vertex corresponding to an $\ell$-link of $Y$. By Lemma \ref{lem_hub}, $H(Y)$ is connected.

By definition each edge of $M$ corresponds to an edge $e$ of $G(\ell)$ between two different branch sets, say $X_i$ and $X_j$. Let $Y$ be the graph consisting of $X_i, X_j$ and $e$. Then $H(X_i)$ and $H(X_j)$ are vertex-disjoint since $X_i$ and $X_j$ are vertex-disjoint. By the analysis above, $H(X_i)$ and $H(X_j)$ are connected subgraphs of the connected graph $H(Y)$. Thus there is a path $Q$ of $H(Y)$ joining $H(X_i)$ and $H(X_j)$ only at end vertices. Further, if $\ell$ is even, then $Q$ is an edge; otherwise, $Q$ is a $2$-path whose middle vertex corresponds to an $\ell$-link $L$ of $Y$ such that $c_L = e$. This implies that $Q$ is internally disjoint with its counterparts and has no inner vertex in any branch set. Then, by contracting each $H(X_i)$ to a vertex, and $Q$ to an edge, we obtain an $M$-minor of $H$.
\qed

Now we are ready to give a lower bound for the Hadwiger number of $\mathbb{L}_\ell(G)$.

\medskip
\noindent{\textbf{Proof of Theorem \ref{thm_hgeqg}}}. Since $H := \mathbb{L}_\ell(G)$ contains an edge, $t := \eta(H) \geqslant 2$.
We first show that $t \geqslant \dege := \dege(G)$. By definition there exists a subgraph $X$ of $G$ of $\delta(X) = \dege$. We may assume that $\dege \geqslant 3$. Then $X$ contains an $(\ell - 1)$-link $P$ such that $\mathscr{L}(P) \neq \emptyset$. By Lemma \ref{lem_VEpartition}, $\mathscr{L}^{[\ell]}(P)$ is the edge set of a complete bipartite subgraph of $H$ with a $K_{\dege - 1, \dege - 1}$-subgraph. By Zelinka \cite{Zelinka1976}, $K_{\dege - 1, \dege - 1}$ contains a $K_{\dege}$-minor. Thus $t \geqslant \dege$ as desired.

We now show that
$t \geqslant \eta := \eta(G)$. If $\eta = 3$, then $G$ contains a cycle $O$ of length at least $3$, and $H$ contains a $K_3$-minor contracted from $\mathbb{L}_\ell(O)$. Now assume that $G$ is connected with $\eta \geqslant 4$. Repeatedly delete vertices of degree $1$ in $G$ until $\delta(G) \geqslant 2$. Then $G = G(\ell)$. Clearly, this process does not reduce the Hadwiger number of $G$. So $G$ contains branch sets of a $K_\eta$-minor covering $V(G)$ (see \cite{Wood2011}). If every branch set contains an $\ell$-link, then the statement follows from Corollary \ref{coro_lminor}. Otherwise, there exists some branch set $X$ with $\diam(X) < \ell$. Since $\eta \geqslant 4$, $Y := G - V(X)$ is connected and contains a cycle. Thus by Lemma \ref{lem_ltConCycle}, $H$ contains a $K_\eta$-minor since $|E(X, Y)| \geqslant \eta - 1$.
\qed



Here we prove Hadwiger's conjecture for $\mathbb{L}_\ell(G)$ for even $\ell \geqslant 2$.

\medskip
\noindent{\textbf{Proof of Theorem \ref{thm_hcllink}(2)}.} Let $\dege := \dege(G)$, $\ell \geqslant 2$ be an even integer, and $H := \mathbb{L}_{\ell}(G)$. By \cite[Proposition 5.2.2]{Diestel2010}, $\chi := \chi(G) \leqslant \dege + 1$. So by Theorem \ref{thm_chrl}, $\chi(H) \leqslant \min\{\dege + 1, \frac{2}{3}\dege + \frac{5}{3}\}$. If $\dege \leqslant 4$, then $\chi(H) \leqslant 5$. By Robertson et al. \cite{RST1993}, Hadwiger's conjecture holds for $H$ in this case. Otherwise, $\dege \geqslant 5$. By Theorem \ref{thm_hgeqg}, $\eta(H) \geqslant \dege \geqslant \frac{2}{3}\dege + \frac{5}{3} \geqslant \chi(H)$ and the statement follows.
\qed

We end this paper by proving Hadwiger's conjecture for $\ell$-link graphs of biconnected graphs for $\ell \geqslant 1$.


\medskip
\noindent{\textbf{Proof of Theorem \ref{thm_hcllink}(1)}.} By Reed and Seymour \cite{Reed2004}, Hadwiger's conjecture holds for $H := \mathbb{L}_\ell(G)$ for $\ell = 1$. By Theorem \ref{thm_hcllink}(2), the conjecture is true if $\ell \geqslant 2$ is even. So we only need to consider the situation that $\ell \geqslant 3$ is odd. If $G$ is a cycle, then $H$ is a cycle and the conjecture holds \cite{Hadwiger1943}. Now let $v$ be a vertex of $G$ with degree $\Delta := \Delta(G) \geqslant 3$. By Theorem \ref{thm_chrl}, $\chi(H) \leqslant \Delta + 1$. Since $G$ is biconnected, $Y := G - v$ is connected. By Lemma \ref{lem_ltConCycle}, if $Y$ contains a cycle, then $\eta(H) \geqslant \Delta + 1 \geqslant \chi(H)$. Now assume that $Y$ is a tree, which implies that $G$ is $K_4$-minor free. By Lemma \ref{lem_ltConnected}, $\eta(H) \geqslant \Delta$. By Theorem \ref{thm_chrl}, $\chi(H) \leqslant \chi' := \chi'(G)$. So it is enough to show that $\chi' = \Delta$.

Let $U := \{u \in V(Y)|\;\deg_Y(u) \leqslant 1\}$. Then $|U| \geqslant \Delta(Y)$. Let $\hat{G}$ be the underlying simple graph of $G$, $t := \deg_{\hat{G}}(v) \geqslant 1$ and $\hat{\Delta} := \Delta(\hat{G}) \geqslant t$. Since $G$ is biconnected, $U \subseteq N_G(v)$. So $t \geqslant |U| \geqslant \Delta(Y)$. Let $u \in U$. When $|U| = 1$, $t = \deg_{\hat{G}}(u) = 1$. When $|U| \geqslant 2$, $\deg_{\hat{G}}(u) = 2 \leqslant |U| \leqslant t$. Thus $t = \hat{\Delta}$. Juvan et al. \cite{JMT1999} proved that the edge-chromatic number of a $K_4$-minor free simple graph equals the maximum degree of this graph. So $\hat{\chi}' := \chi'(\hat{G}) = \hat{\Delta}$ since $\hat{G}$ is simple and $K_4$-minor free. Note that all parallel edges of $G$ are incident to $v$. So $\chi' = \hat{\chi}' + \deg_G(v) - t = \hat{\Delta} + \Delta - \hat{\Delta} = \Delta$ as desired.
\qed



%% file: 41Introduction.tex
\section{Introduction and main results}
As a generalisation of line graphs \cite{Whitney1932} and path graphs \cite{BH1989}, the $\ell$-link graph of a given graph was introduced by Jia and Wood \cite{JiaWood2013} who studied the connectedness, chromatic number and minors of $\ell$-link graphs based on the structure of the given graph. This paper deals with the reverse; that is, for an integer $\ell \geqslant 0$ and a given finite graph $H$, we study the graphs whose $\ell$-link graphs are isomorphic to $H$.

Unless stated otherwise, all graphs are undirected and loopless. A graph may be finite or infinite, and may be simple or contain parallel edges. In particular, $H$ always denotes a finite graph. The \emph{order} and \emph{size} of $H$ are \emph{$n(H) := |V(H)|$} and \emph{$m(H) := |E(H)|$} respectively. Throughout this paper, $\ell \geqslant 0$ is an integer. An \emph{$\ell$-link} is a walk of length $\ell$ in which consecutive edges are different. We identify an $\ell$-link with its reverse sequence. In particular, an \emph{$\ell$-path} is an $\ell$-link without repeated vertices. The \emph{$\ell$-link graph $\mathbb{L}_\ell(G)$} of a graph $G$ is defined to have vertices the $\ell$-links of $G$, and two vertices are adjacent if their corresponding $\ell$-links form an $(\ell + 1)$-link of $G$. If further $G$ contains parallel edges, then two $\ell$-links may form $\mu \geqslant 2$ different $(\ell + 1)$-links \cite{JiaWood2013}. In this case, we give $\mu$ edges between the two corresponding vertices in $\mathbb{L}_\ell(G)$. More strict definitions can be found in Section \ref{sec_term}.

A graph $G$ is an \emph{$\ell$-root} of $H$ if $\mathbb{L}_\ell(G) \cong H$. Let $\mathbb{R}_\ell(H)$ be the set of minimal (up to subgraph relation) $\ell$-roots of $H$. Before presenting the main results, we would like to mention some work about the number of minimal $\ell$-roots. By definition $\mathbb{R}_0(H) = \{H\}$. Lemma \ref{lem_barK2} indicates that $|\mathbb{R}_\ell(\bar{K}_2)|$ increases with $\ell$. We prove in Lemma \ref{lem_copylemma} that, for a fixed $\ell \geqslant 4$ and any given number $k$, there exists a graph $H$ with $|\mathbb{R}_\ell(H)| > k$. The \emph{line graph $\mathbb{L}(G)$} of $G$ is the simple graph with vertex set $E(G)$, in which two vertices are adjacent if their corresponding edges share a common end vertex in $G$. By \cite[Remark 2.1]{JiaWood2013}, $\mathbb{L}_1(G) = \mathbb{L}(G)$ if and only if $G$ is simple. So Whitney's theorem \cite{Whitney1932} can be restated as: for a simple connected finite graph $H$, $\mathbb{R}_1(H) = \{K_3, K_{1,3}\}$ if $H \cong K_3$, and $|\mathbb{R}_1(H)| \leqslant 1$ otherwise. The theorem below is a qualitative generalisation of Whitney's theorem.

\begin{thm}\label{thm_minirootbounded} Let $\ell \geqslant 0$ be an integer, and $H$ be a finite graph. Then
the maximum degree, order, size, and total number of minimal $\ell$-roots of $H$ are finite and bounded by functions of $H$ and $\ell$.
\end{thm}

The notions of tree-decomposition and tree-width were studied extensively by Robertson and Seymour in proving that finite graphs are well-quasi-ordered by the minor relation \cite{GMXX}. Let $G$ be a hypergraph, $T$ be a tree, and $\mathcal{V} := \{V_w | w \in V(T)\}$ be a set cover of $V(G)$. The pair $(T, \mathcal{V})$ is called a \emph{tree-decomposition} of $G$ if, first of all, for every $e \in E(G)$, there exists $V \in \mathcal{V}$ containing all vertices incident to $e$. And secondly, for every path $[w_0, \ldots, w_i, \ldots, w_\ell]$ of $T$, we have $V_{w_0} \cap V_{w_\ell} \subseteq V_{w_i}$. The \emph{width} of the tree decomposition $(T, \mathcal{V})$ is \emph{$\tw(T, \mathcal{V}) := \sup\{|V| - 1| V \in \mathcal{V}\}$}. The \emph{tree-width $\tw(G)$} is the minimum width over the tree-decompositions of $G$. The \emph{tree-diameter $\tdi(G)$} is the minimum diameter of $T$ over the tree-decompositions $(T, \mathcal{V})$ of $G$ with width $\tw(G)$.

\begin{thm}\label{thm_llinkbqo} Let $\ell \geqslant 0$ be an integer, and $H$ be a finite graph. Then
the tree-width and tree-diameter of the $\ell$-roots of $H$ are finite and bounded by functions of $H$ and $\ell$. Further, the $\ell$-roots of a finite graph are better-quasi-ordered by the induced subgraph relation.
\end{thm}


Apply Theorem \ref{thm_llinkbqo} to the $\ell$-roots of the \emph{null graph}, we have:
\begin{lem}\label{lem_treeBetterInsub}
The trees of bounded diameter are better-quasi-ordered by the induced subgraph relation.
\end{lem}








We use \emph{$\mathbb{R}_\ell[H]$} to denote the set of $\ell$-roots of $H$. To see an example, let \emph{$M_G(u, v)$} be the set of edges of $G$ between $u, v \in V(G)$. Then for each $e \in M_G(u, v)$, the $1$-link $[u, e, v]$ is formed by the $0$-links $[u]$ and $[v]$ in $G$. Hence $\mathbb{L}_0(G) \cong G$ and $\mathbb{R}_0[G] = \{G\}$.

Denote by \emph{$X \subseteq Y$}, \emph{$X \subset Y$}, \emph{$X \leqslant Y$}, and \emph{$X < Y$} that $X$ is isomorphic to a subgraph, proper subgraph, induced subgraph and proper induced subgraph of a graph $Y$. A graph is \emph{$\ell$-finite} if its $\ell$-link graph is finite. So all finite graphs are $\ell$-finite, but not vice versa. For example, let $T^{t}$ be the tree obtained by pasting the middle vertex of a $4$-path at the center of a star $K_{1, t}$. Then $T^{\infty}$ is infinite, and is $4$-finite since its $4$-link graph is $K_1$.

Two $\ell$-finite graphs
$X$ and $Y$ are \emph{$\ell$-equivalent}, written \emph{$X \thicksim_\ell Y$}, if there exists a graph $Z \subseteq X, Y$ such that $\mathbb{L}_\ell(X) \cong \mathbb{L}_\ell(Y) \cong \mathbb{L}_\ell(Z)$. For every pair of integers $i, j \geqslant 0$, we have $T^i \thicksim_4 T^j$ since $\mathbb{L}_4(T^i) \cong \mathbb{L}_4(T^j) \cong \mathbb{L}_4(T^0) \cong K_1$.
An $\ell$-finite graph $X$ is \emph{$\ell$-minimal} if $X$ is null or $\mathbb{L}_\ell(Y) \subset \mathbb{L}_\ell(X)$ for every $Y \subset X$. For instance, an $\ell$-path is $\ell$-minimal. By definition a graph is $\ell$-minimal if and only if it is a minimal $\ell$-root of a finite graph.

The following two lemmas are proved in Section \ref{sec_consEquiClass}. Let $H$ be a finite graph. Lemma \ref{lem_miniUniq} states that $\mathbb{L}_\ell[H]$ is the union of

\begin{lem}\label{lem_miniUniq}
For each integer $\ell \geqslant 0$, $\sim_\ell$ is an equivalence relation on $\ell$-finite graphs, such that each $\ell$-equivalence class contains a unique (up to the isomorphism) $\ell$-minimal graph. And this graph is isomorphic to an induced subgraph of every graph in its class.
\end{lem}

Our next lemma tells that an $\ell$-root of a finite graph $H$ is a certain combination of a minimal $\ell$-root and trees of bounded diameter. This transfers the study of $\ell$-roots into that of minimal $\ell$-roots. The minimal $\ell$-roots of a cycle are obtained in Section \ref{sec_example}, where all $\ell$-roots of a cycle are also characterised by applying Lemma \ref{lem_conClass}.


\begin{lem}\label{lem_conClass}
Let $\ell \geqslant 0$ be an integer, and $G$ be the minimal graph of an $\ell$-equivalence class. Then a graph belongs to this class if and only if it can be obtained from $G$ by:
\begin{itemize}
\item[{\bf (1)}] For each acyclic component $T$ of $G$ of diameter at most $2\ell - 2$, and every vertex $u$ of eccentricity $s \leqslant \ell - 2$ in $T$, pasting to $u$ the root of a rooted tree of height at most $\ell - s - 1$.
\item[{\bf (2)}] Adding to $G$ zero or more acyclic components of diameter at most $\ell - 1$.
\end{itemize}
\end{lem}





Introduced by Broersma and Hoede \cite{BH1989}, the \emph{$\ell$-path graph $\mathbb{P}_\ell(G)$} is the simple graph with vertices the $\ell$-paths of $G$, where two vertices are adjacent if the union of their corresponding paths forms a path or a cycle of length $\ell + 1$ in $G$. By Jia and Wood \cite{JiaWood2013}, when $\ell \geqslant 2$, we have $\mathbb{P}_\ell(G) \leqslant \mathbb{L}_\ell(G)$, where the equation holds if and only if $\girth(G) > \ell$. We say $G$ is an \emph{$\ell$-path root} of $H$ if $\mathbb{P}_\ell(G) \cong H$. Let $\mathbb{Q}_\ell(H)$ be the set of minimal (up to subgraph relation) $\ell$-path roots of $H$. Li \cite{LiXL1996} proved that $H$ has at most one simple $2$-path root of minimum degree at least $3$. Prisner \cite{Prisner2000} showed that $\mathbb{Q}_\ell(H)$ contains at most one simple graph of minimum degree greater than $\ell$.
By Li and Liu \cite{LiLiu2008}, if $H$ is connected and nonnull, then $\mathbb{Q}_2(H)$ contains at most two simple graphs. In fact the finite graphs having exactly two simple minimal $2$-path roots have been characterised by Aldred, Ellingham, Hemminger and Jipsen \cite{AEHJ1997}. Some results about $\ell$-roots can be proved, with slight variations, for $\ell$-path roots:

\begin{thm}\label{thm_pathrootsBound} Let $\ell \geqslant 0$ be an integer, and $H$ be a finite graph. Then
the order, size, and total number of minimal $\ell$-path roots of $H$ are finite and bounded by functions of $H$ and $\ell$. Further, the tree-width and tree-diameter of $\ell$-path roots of $H$ are finite and bounded by functions of $H$ and $\ell$. Moreover, the $\ell$-path roots (respectively, of bounded multiplicity) of a finite graph are better-quasi-ordered by the (respectively, induced) subgraph relation.
\end{thm}

Ding \cite{Ding1992} proved that finite simple $\ell$-path-free graphs are well-quasi-ordered by the induced subgraph relation. Apply Theorem \ref{thm_pathrootsBound} to the $\ell$-path roots of the null graph, we generalise Ding's theorem as follows:
\begin{lem}\label{lem_lpathfreebqo}
Given a finite graph $H$, the $H$-minor free graphs with or without loops (respectively, of bounded multiplicity) are better-quasi-ordered by the (respectively, induced) subgraph relation if and only if $H$ is a disjoint union of paths.
\end{lem}

%% file: 42Terminology.tex
\section{Terminology}\label{sec_term}
We list in this section some necessary definitions and simple facts.
Let $G$ be a graph, and \emph{$c(G)$} (respectively, \emph{$o(G)$}, \emph{$a(G)$}) be the cardinality of the set of (respectively, cyclic, acyclic) components of $G$.
The \emph{distance $\dist_G(u, v)$} between $u, v \in V(G)$ is $+ \infty$ if $u, v$ are in different components of $G$, and the minimum length of a path of $G$ between $u, v$ otherwise. The \emph{eccentricity} of $v \in V(G)$ is \emph{$\ecc_G(v) := \sup\{\dist_G(u, v) | u \in V(G)\}$}. The \emph{diameter} of $G$ is \emph{$\diam(G) := \sup\{\ecc_G(v)| v \in V(G)\}$}. The \emph{radius} of $G$ is \emph{$\radi(G) := \min\{\ecc_G(v)| v \in V(G)\}$}. Clearly, for each tree $T$, we have $\radi(T) = \lceil \diam(T)/2\rceil$. Denote by \emph{$K_t$} the complete graph on $t \geqslant 0$ vertices. In particular, \emph{$K_0$} is called the \emph{null} graph. Denote by $tG$ is the disjoint union of $t \geqslant 0$ copies of $G$. For $t \geqslant 1$, \emph{$\bar{K}_t := tK_1$} is called the \emph{empty graph} on $t$ vertices. For $s \geqslant 1$, the \emph{$s$-subdivision $G^{\langle s\rangle}$} of $G$ is the graph obtained by replacing every edge of $G$ with an $s$-path. So $G^{\langle 1\rangle} = G$. Let $e$ be an edge of a tree $T$ with end vertices $u$ and $v$. Let \emph{$T_e^u$} be the component of $T - e$ containing $u$, and \emph{$T_u^e := T_e^v \cup \{e\}$}. A \emph{unit} is a vertex or an edge. The subgraph of $G$ induced by $V \subseteq V(G)$ is the maximal subgraph of $G$ with vertex set $V$. For $\emptyset \neq E \subseteq E(G)$, the subgraph of $G$ induced by $E \cup V$ is the minimal subgraph of $G$ with edge set $E$, and vertex set including $V$.

An \emph{$\ell$-arc} (or \emph{$*$-arc} if we ignore the length) is a sequence $\vec{L} := (v_0, e_1, \ldots, e_\ell, v_\ell)$, where $e_i$ is an edge of end vertices $v_{i - 1}$ and $v_i$ such that $e_{j} \neq e_{j + 1}$ for $i \in [\ell] := \{1, 2, \ldots, \ell\}$ and $j \in [\ell - 1]$. Note that $\vec{L}$ is different from $-\vec{L} := (v_{\ell}, e_\ell, \ldots, e_1, v_0)$ unless $\ell = 0$. For each $i \in [\ell]$, $\vec{e}_i := (v_{i - 1}, e_i, v_i)$ is called an \emph{arc} for short. $v_0$, $v_\ell$, $\vec{e}_1$ and $\vec{e}_\ell$ are \emph{tail vertex}, \emph{head vertex}, \emph{tail arc} and \emph{head arc} of $\vec{L}$ respectively.   The $\ell$-link (or $*$-link if we ignore the length) $L := [v_0, e_1, \ldots, e_\ell, v_\ell] = [v_{\ell}, e_\ell, \ldots, e_1, v_0]$ is obtained by taking $\vec{L}$ and $-\vec{L}$ as a single object; that is, $L := \{\vec{L}, -\vec{L}\}$. For $0 \leqslant i \leqslant j \leqslant \ell$, $\vec{R} := \vec{L}(i, j) := (v_i, e_{i + 1}, \ldots, e_j, v_j)$ is called an \emph{$(j - i)$-arc} (or a \emph{subsequence} for short) of $\vec{L}$, and $\vec{L}[i, j] := R$ is an \emph{$(j - i)$-link} (or a \emph{subsequence} for short) of $L$. For example, a $0$-link is a vertex, and a $1$-link can be identified with an edge. For $\ell \geqslant 1$, we say $L$ is \emph{formed} by the $(\ell - 1)$-links $\vec{L}[0, \ell - 1]$ and $\vec{L}[1, \ell]$. An \emph{$\ell$-dipath} is an $\ell$-arc without repeated vertices. We say $\vec{L}$ is an \emph{$\ell$-dicycle} if $v_0 = v_\ell$ and $\vec{L}(0, \ell - 1)$ is a dipath. An \emph{$\ell$-path} is an $\ell$-link without repeated vertices. We use \emph{$\vec{\mathscr{L}}_\ell(G)$}, \emph{$\mathscr{L}_\ell(G)$}, and \emph{$\mathscr{P}_\ell(G)$} to denote the sets of $\ell$-arcs, $\ell$-links, and $\ell$-paths of $G$ respectively.

Let $\vec{L} := (v_0, e_1, \ldots, e_\ell, v_\ell) \in \vec{\mathscr{L}}_\ell(G)$, and $\vec{R} := (u_0, f_1, \ldots, f_s, u_s) \in \vec{\mathscr{L}}_s(G)$ such that $v_\ell = u_0$ and $e_\ell \neq f_1$. The \emph{conjunctions} of $\vec{L}$ and $\vec{R}$ are $\vec{Q} := (\vec{L}.\vec{R}) := (v_0, e_1, \ldots, e_\ell, v_\ell = u_0, f_1, \ldots, f_s, u_s) \in \vec{\mathscr{L}}_{\ell + s}(G)$ and $[\vec{L}.\vec{R}] := Q \in \mathscr{L}_{\ell + s}(G)$. For $\vec{Q} \in \vec{\mathscr{L}}_{\ell + s}(G)$, let $\vec{L}_i := \vec{Q}(i, \ell + i)$, and $\vec{Q}_j := \vec{Q}(j - 1, \ell + j)$ for $i \in \{0, 1, \ldots, s\}$ and $j \in [s]$. By definition $Q_j \in \mathscr{L}_{\ell + 1}(G)$ yields an edge $Q_j^{[\ell]} := [L_{j - 1}, Q_j, L_j]$ of $\mathbb{L}_\ell(G)$. So $Q^{[\ell]} := [L_0, Q_1, L_1, \ldots, L_{s - 1}, Q_s, L_s]$ can be seen as an $s$-link, while $\vec{Q}^{[\ell]} := (L_0, Q_1, L_1, \ldots, L_{s - 1}, Q_s, L_s)$ is an $s$-arc of $\mathbb{L}_\ell(G)$. We say that $L_0$ can be shunted to $L_s$ through $\vec{Q}$. $Q^{\{\ell\}} := \{L_0, L_1, \ldots, L_s\}$ and $\vec{Q}^{\{\ell\}} := \{\vec{L}_0, \vec{L}_1, \ldots, \vec{L}_s\}$ are the sets of \emph{images} of $L_0$ and $\vec{L}_0$ respectively during this shunting. More generally, for $R, R' \in \mathscr{L}_\ell(G)$, we say $R$ can be \emph{shunted} to $R'$ if there are $\ell$-links $R = R_0, R_1, \ldots, R_s = R'$, and $*$-arcs $\vec{P}_1, \ldots, \vec{P}_s$ of $G$ such that $R_{i - 1}$ can be shunted to $R_i$ through $\vec{P}_i$ for $i \in [s]$.

Let $(T, \mathcal{V})$ be a tree-decomposition of $G$, $[v_0, e_1, \ldots, v_\ell] \in \mathscr{P}_\ell(T)$, and \emph{$\mathcal{V}_T(v_0, v_\ell)$} be the set of minimal sets, up to subset relation, among $V_{v_i}$ and $V_{e_j} := V_{v_{j - 1}} \cap V_{v_j}$ for $i \in \{0, \ldots, \ell\}$ and $j \in [\ell]$. $(T, \mathcal{V})$ is a \emph{linked tree-decomposition} if for every pair of $u, v \in V(T)$, and $U \subseteq V_u$, $V \subseteq V_v$ such that $|U| = |V| =: k$, either $G$ contains $k$ disjoint paths from $U$ to $V$, or there exists some $W \in \mathcal{V}_T(u, v)$ such that $|W| < k$. Kruskal's theorem \cite{Kruskal1960} states that finite trees are well-quasi-ordered by the topological minor relation. Nash-Williams \cite{NashWilliams1965} generalised this theorem and proved that infinite trees are \emph{better-quasi-ordered} under the same relation. Let $\mathcal{A}$ be the set of all finite ascending sequences of nonnegative integers. For $A, B \in \mathcal{A}$, written $A <_{\mathcal{A}} B$ if $A$ is a strict initial subsequence of some $C \in \mathcal{A}$, and by deleting the first term of $C$, we obtain $B$. Let $\mathcal{B}$ be an infinite subset of $\mathcal{A}$, and $\bigcup\mathcal{B}$ be the set of nonnegative integers appeared in some sequence of $\mathcal{B}$. $\mathcal{B}$ is called a block if it contains an initial subsequence of every infinite increasing sequence of $\bigcup\mathcal{B}$. Let $\mathcal{Q}$ be a set with a quasi-ordering $\leqslant_{\mathcal{Q}}$. A \emph{$\mathcal{Q}$-pattern} is a function from a block $\mathcal{B}$ into $\mathcal{Q}$. A $\mathcal{Q}$-pattern $\varphi$ is \emph{good} if there exist $A, B \in \mathcal{B} \subseteq \mathcal{A}$ such that $A <_{\mathcal{A}} B$ and $\varphi(A) \leqslant_{\mathcal{Q}} \varphi(B)$. $\mathcal{Q}$ is said to be \emph{better-quasi-ordered} by $\leqslant_{\mathcal{Q}}$ if every $\mathcal{Q}$-pattern is good. For $j \geqslant 1$, define a quasi-ordering on $\mathcal{Q}^j$ as: $(q_1, \ldots, q_j) \leqslant_{\mathcal{Q}^j} (q'_1, \ldots, q'_j)$ if $q_i \leqslant_{\mathcal{Q}} q'_i$ for every $i \in [j]$. The following lemma follows directly from Galvin-Prikry theorem  \cite{GalvinPriskry1973}.
\begin{lem}\label{lem_coverbqo}
Let $k \geqslant 1$ be an integer. Then
$\mathcal{Q} = \bigcup_{i = 1}^k\mathcal{Q}_i$ is better-quasi-ordered if and only if for every $i \in [k]$, $\mathcal{Q}_i$ is better-quasi-ordered if and only if for every $j \geqslant 1$, $\mathcal{Q}^j$ is better-quasi-ordered.
\end{lem}



Define a quasi-ordering on a set $\mathcal{S}$ of sequences of $\mathcal{Q}$: $S_1 \leqslant_\mathcal{S} S_2$ if there is an order-preserving injection $\varphi: S_1 \mapsto S_2$ such that $q \leqslant \varphi(q)$ for every $q \in S_1$. The lemma below is due to Nash-Williams \cite{NashWilliams1968}.

\begin{lem}\label{lem_subseqbqo}
Any finite quasi-ordered set is better-quasi-ordered.
$\mathcal{Q}$ is better-quasi-ordered if and only if any set of sequences on $\mathcal{Q}$ is better-quasi-ordered.
\end{lem}

%


%% file: 43Examples.tex
\section{Examples and basis}\label{sec_example}
We begin with some examples and basic analysis which help to build some general impressions on $\ell$-roots, and explain some of our motivations.


First of all, we characterise the minimal $\ell$-roots of $\bar{K}_2$.
\begin{lem}\label{lem_barK2}
Let $P := [v_0, \ldots, v_\ell]$ be an $\ell$-path, and $T_i$ be obtained from $P$ by pasting $v_i$ at an end vertex of another $i$-path. Then $\mathbb{R}_\ell(\bar{K}_2) = \{2P, T_i | 1 \leqslant i \leqslant \lfloor\frac{\ell - 1}{2}\rfloor\}$. Further, $|\mathbb{R}_\ell(\bar{K}_2)|$ is $1$ if $\ell = 0$, and is $\lfloor\frac{\ell + 1}{2}\rfloor$ if $\ell \geqslant 1$.
\end{lem}
\pf Clearly, for $1 \leqslant i \leqslant \lfloor\frac{\ell - 1}{2}\rfloor$, $\mathbb{L}_\ell(2P) \cong \mathbb{L}_\ell(T_i) \cong \bar{K}_2$. If $G \in \mathbb{R}_\ell(\bar{K}_2)$ contains a cycle $O$, then $\mathbb{L}_\ell(G)$ contains a cycle $\mathbb{L}_\ell(O)$, which is impossible. Thus $G$ is a forest containing exactly one $\ell$-path $Q$ other than $P$. If $P$ and $Q$ are vertex disjoint, then $G = P \cup Q \cong 2P$ because of the minimality. Otherwise, assign directions such that $\vec{P} = (\vec{P}_1.\vec{R}.\vec{P}_2)$ and $\vec{Q} = (\vec{Q}_1.\vec{R}.\vec{Q}_2)$, where $R$ is a maximal common path of $P$ and $Q$, $P_i \in \mathscr{P}_{s_i}(G)$ and $Q_i \in \mathscr{P}_{t_i}(G)$ for $i \in \{1, 2\}$. Since $P \neq Q$, without loss of generality, $\vec{P}_1 \neq \vec{Q}_1$ and $s_1 \geqslant t_1$. Then $s_2 \leqslant t_2$, and $\vec{L} := (\vec{P}_1.\vec{R}.\vec{Q}_2) \in \vec{\mathscr{L}}_{\ell + t_2 - s_2}(G)\setminus\{\vec{Q}\}$. Since $\bar{K}_2$ contains no edge, $\mathscr{L}_{\ell + 1}(G) = \emptyset$ and so $s_2 = t_2$. Thus $s_1 = t_1 \geqslant 1$, and $L \in \mathscr{P}_\ell(G)$. So $\vec{L} = \vec{P}$ since otherwise, $G$ contains three pairwise different $\ell$-paths $L, P$ and $Q$. Note that $[\vec{P}_1.-\vec{Q}_1] \in \mathscr{P}_{2s_1}(G)\setminus \{P, Q\}$. So $2s_1 < \ell$ and the lemma follows.
\qed

%


By Whitney \cite{Whitney1932}, $\mathbb{R}_1(K_3) = \{K_3, K_{1, 3}\}$. As a generalisation, Broersma and Hoede \cite{BH1989} pointed out that a $6$-cycle is the $2$-path (and hence $2$-link) graph of itself and $K_{1, 3}^{\langle 2\rangle}$. Below we characterise the minimal $\ell$-roots of all cycles. Clearly, for a given $\ell \geqslant 0$, every cycle has a unique cyclic minimal $\ell$-root which is isomorphic to itself. So we only need to consider acyclic minimal $\ell$-roots.




\begin{lem}\label{lem_cyclerootT}
Let $T$ be a minimal acyclic $\ell$-root of a $t$-cycle. Then $\ell \geqslant 1$, and either $t = 3\ell$ and $T \cong K_{1, 3}^{\langle\ell\rangle}$, or there is $s \geqslant 1$ such that $t = 4s$, $\ell \geqslant 2s + 1$, and $T$ is obtained by joining the middle vertices of two $2s$-paths by an $(\ell - s)$-path.
\end{lem}
\pf Since $O := \mathbb{L}_\ell(T)$ is $2$-regular, $D(u, v) := \deg_T(u) + \deg_T(v) = 4$ if $\dist_T(u, v) = \ell$, and $D(u, v) \leqslant 4$ if $\dist_T(u, v) \geqslant \ell$. By the minimality of $T$, for each leaf $w$ of $T$, there exists some $v \in V(T)$ such that $\dist_T(v, w) = \ell$ and $\deg_T(v) = 3$.
By Lemma \ref{lem_Tdc}, for every $v \in V(T)$ with $\ecc_T(v) < \ell$, $\deg_T(v) \leqslant c(O) + 1 \leqslant 2$. So $\deg_T(v) \in \{1, 2, 3\}$ for $v \in V(T)$, and $T$ contains $k \geqslant 1$ vertices of degree $3$. If $k = 1$, then $T$ contains exactly three leaves (\cite[Page 67]{ChartrandOellermann1993}), each has distance $\ell$ with $v$. So $T \cong K_{1, 3}^{\langle\ell\rangle}$ and $t = 3\ell$. If $k \geqslant 2$, there exists $q \in [\ell - 1]$ and $\vec{Q} := (v_0, \ldots, v_q) \in \vec{\mathscr{L}}_q(T)$ such that $\deg_T(v_i) = 3$ for $i \in \{0, q\}$. If $k = 2$, there are four leaves \cite[Page 67]{ChartrandOellermann1993} of distance $\ell$ with $v_0$ or $v_q$. So $T$ is the union of two paths $[\vec{L}_i.\vec{Q}.\vec{R}_i]$, where $i \in \{1, 2\}$, and $L_i, R_i$ are four internally disjoint paths of length $\ell_i, s_i \geqslant 1$ respectively. Consider $\vec{L}_i := (w_0, \ldots, w_{\ell_i} = v_0)$. If $\dist_T(w_0, v_0) = \ell$, then $\dist_T(w_1, v_q) = \ell - 1 + q \geqslant \ell$, and $D(w_1, v_q) > 4$, a contradiction. So by the analysis above, $\dist_T(w_0, v_q) = \ell$. Thus $s := \ell_i = s_i \geqslant 1$ for $i \in \{1, 2\}$ and $q = \ell - s$. Moreover, let $\vec{P}_{ij} := (\vec{L}_i.\vec{Q}.\vec{R}_j)$ for $i, j \in \{1, 2\}$. Then the images during the shunting of $L$ through $\vec{P}_{11}, -\vec{P}_{21}, \vec{P}_{22}, -\vec{P}_{12}$ form a $4s$-cycle.
Note that $[\vec{L}_1.\vec{L}_2] \in \mathscr{P}_{2s}(T)$ is not an image mentioned above. So $2s < \ell$ and the lemma follows in this case. We still need to show that $k < 3$. Otherwise, there exists some $p \in [q - 1]$ of $\deg_T(v_p) = 3$. By the analysis above, $\ecc_T(v_p) \geqslant \ell$, and so there exists $\vec{L} := (u_0, \ldots, u_\ell = v_p) \in \vec{\mathscr{L}}_\ell(T)$ such that $v_0, v_q$ are separated from $u_0$ by $v_p$. Then $\dist_T(u_1, v_q) = \ell - 1 + q - p \geqslant \ell$, and $D(u_1, v_q) > 4$, a contradiction.
\qed
\noindent {\bf Remark.} Lemma \ref{lem_conClass} and \ref{lem_cyclerootT} provide us all $\ell$-roots $G$ of a $t$-cycle: if $G$ is cyclic, it is the disjoint union of a $t$-cycle and zero or more trees of diameter at most $\ell - 1$. Let $G$ be a forest and $\ell \geqslant 1$. If $t = 3\ell$, since $\diam(K_{1, 3}^{\langle\ell\rangle}) = 2\ell$, $G$ is the disjoint union of $K_{1, 3}^{\langle\ell\rangle}$ and trees of diameter at most $\ell - 1$. In the final case of Lemma \ref{lem_cyclerootT}, $\diam(T) = \ell + s \leqslant 2\ell - 2$. Let $[v_0, \ldots, v_{\ell - s}]$ be the path of $T$ between the middle vertices of the two $2s$-paths. Then $\ecc_{T}(v_i) = \max\{i, \ell - i\}$ for $i \in [\ell - s]$. So $G$ is obtained from $T$ by first pasting to each $v_i$, where $i \in \{2, 3, \ldots, \ell - s - 2\}$, the root of a rooted tree of height less than $\min\{i, \ell - s - i\}$, and then adding acyclic components of diameter less than $\ell$.

From Lemma \ref{lem_cyclerootT} and our next example, a $4s$-cycle has at least three minimal $(2s + 1)$-path roots, two of which are cyclic, where $s \geqslant 1$ is an integer.
\begin{exa}
Let $s \geqslant 1$ and $\ell \geqslant s + 1$ be integers. Let $G(s, \ell)$ be the graph formed by connecting two $(s + 1)$-cycles with an $(\ell - s)$-path. One can easily check that $G(s, \ell)$ is a minimal $\ell$-path root of a $4s$-cycle.
\end{exa}

Broersma and Hoede \cite{BH1989} asked that, for $\ell = 2$, whether there exist three pairwise non-isomorphic simple connected graphs whose $\ell$-path graphs are isomorphic to the same connected nonnull graph. A negative answer was given by Li and Liu \cite{LiLiu2008}. The following lemma addresses this problem in the case of $\ell \geqslant 3$. It gives a positive answer by showing that, there exist infinite many trees $T$ of diameter $3$ such that $\mathbb{Q}_\ell(T) \cap \mathbb{R}_\ell(T)$ contains at least four trees. Further, let $\ell \geqslant s \geqslant 4$ and $k \geqslant 0$ be given integers. Then there exists a tree $T$ of diameter $s$, such that $\mathbb{Q}_\ell(T) \cap \mathbb{R}_\ell(T)$ contains at least $k$ trees.

\begin{lem}\label{lem_copylemma}
Let $T$ be a finite tree and $v \in V(T)$ of degree $d$. Assign to $v$ an integer $t_v$ as: if $d \geqslant 2$, $t_v := \diam(T)$. If $d \leqslant 1$, on one hand, if $T$ is a path, then $t_v := -1$. On other hand, there exists a path $[v, \ldots, e, u]$ of minimum length such that $\deg_T(u) \geqslant 3$. In this case, let $t_v := \diam(T_e^u)$. Denote by $T(v, \ell)$ the tree obtained by pasting an end vertex of an extra $\ell$-path to $v$. Then for each $\ell \geqslant t_v + 1$, $T \cong \mathbb{L}_\ell(T(v, \ell)) = \mathbb{P}_\ell(T(v, \ell))$.
\end{lem}
\pf Let $\vec{L}$ be the $\ell$-arc of head vertex $v$ such that $L$ is the extra path. Consider the shunting of $L$ in $G := T(v, \ell)$. One can check that the mapping $\vec{L}'[\ell, \ell] \mapsto L'$, for every image $\vec{L}'$ of $\vec{L}$, is an isomorphism from $T$ to $\mathbb{L}_\ell(G)$.
\qed
\noindent {\bf Remark.} For a fixed $\ell \geqslant 1$, the number of non-isomorphic $T(v, \ell)$, over all $v \in V(T)$, equals the number of \emph{orbits} of $V(T)$ under the \emph{automorphism group $Aut(T)$} (see \cite{Biggs1993}) of $T$, which is at least $\lfloor \diam(T)/2\rfloor + 1$, with equation holds if and only if the set of leaves is $Aut(T)$-transitive. Further, if $\diam(T) = 3$, then $V(T)$ has up to four $Aut(T)$-orbits. For each $s \geqslant 4$ and $k \geqslant 1$, let $T^k$ be obtained by pasting a leaf of each star $K_{1, i}$, where $i \in [k]$, at the same end vertex of an $(s - 2)$-path. Then $\diam(T) = s$, and the number of $Aut(T^k)$-orbits of $V(T^k)$ is $\lfloor\frac{s}{2}\rfloor + 1$ if $k = 1$, and $s + 2k - 1$ if $k \geqslant 2$.

%% file: 44IncidentL-links.tex
\section{Constructing $\ell$-equivalence classes}\label{sec_consEquiClass}
In this section, we explain the process of constructing $\ell$-roots from minimal $\ell$-roots, which allows us to concentrate on the latter in our future study.

\subsection{Incidence units}
Two $*$-links of a graph $G$ are \emph{incident} if one is a subsequence of the other. A $*$-link is said to be \emph{$\ell$-incident} if it is incident to an $\ell$-link. It follows from the definitions immediately that every $\ell$-link is $\ell$-incident, and every $\ell$-incident $*$-link is $s$-incident, for $s \leqslant \ell$. Conversely, a $t$-link is not $\ell$-incident if and only if it is not $s$-incident for any $s \geqslant \ell$. And if this is the case, $\ell \geqslant t + 1$. In Lemma \ref{lem_barK2}, all units of $T_i$ are $\ell$-incident. However, when $\ell \geqslant 3$, $[\vec{P}_1.-\vec{Q}_1]$ is a $2i$-path that is not $\ell$-incident in $T_i$. A \emph{ray} is an infinite graph with vertex set $\{v_0, v_1, \ldots\}$ and edges $e_i$ between $v_{i - 1}$ and $v_i$, for $i \geqslant 1$. The fact below allows us to focus on incidence units of trees of finite diameters.


\begin{lem}\label{lem_cycleIncident}
Let $G$ be a connected nonnull graph. Then $G$ contains a cycle or a ray if and only if for every $\ell \geqslant 0$, all units of $G$ are $\ell$-incident.
\end{lem}
\pf
$(\Leftarrow)$ Suppose not. Then $G$ is a tree of finite diameter $s$. Then no unit of $G$ is $(s + 1)$-incident. $(\Rightarrow)$ Let $X$ be a cycle or a ray in $G$. Clearly, every unit of $X$ is $\ell$-incident. So we only need to show that every $e \in E(G)\setminus E(X)$ is $\ell$-incident. Since $G$ is connected, there exists a dipath $\vec{P}$ of minimum length with tail edge $e$ and head vertex $x \in V(X)$. Clearly, $X$ contains an $\ell$-arc $\vec{R}$ starting from $x$. Then $L := (\vec{P}.\vec{R})[0, \ell]$ is an $\ell$-link of $G$ incident to $e$.
\qed


The following simple fact tells that, to study $\ell$-incident units of a tree $T$ of finite, we can assume that $\diam(T) \geqslant \ell \geqslant \max\{4, \radi(T) + 1\}$.
\begin{obs}\label{obs_dlr4}
Let $T$ be a tree of finite diameter. For each $\ell \leqslant \radi(T)$ or $\ell \leqslant \min\{3, \diam(T)\}$, all units of $T$ are $\ell$-incident.
\end{obs}



Wu et al. \cite{WuChao2004} presented a linear time algorithm computing the eccentricity of a vertex of a finite tree. Based on this work, the following observation provides a linear time algorithm testing if a vertex is $\ell$-incident in a finite tree.
\begin{obs}\label{obs_vertexincidence}
Let $T$ be a tree and $\ell \geqslant 0$. Then $u \in V(T)$ is $\ell$-incident in $T$ if and only if either $u$ is a leaf and $\ecc_T(u) \geqslant \ell$, or there exist different $e, f \in E(T)$ incident to $u$, such that $\ecc_{T_u^e}(u) + \ecc_{T_u^f}(u) \geqslant \ell$.
\end{obs}

Based on Observation \ref{obs_vertexincidence}, the lemma below can be formalised into a linear time algorithm for testing if an edge of a finite tree is $\ell$-incident.
\begin{lem}\label{lem_elcovered}
Let $\ell \geqslant 0$, and $P$ be a path of a tree $T$. Then all units of $P$ are $\ell$-incident in $T$ if and only if both ends of $P$ are $\ell$-incident in $T$.
\end{lem}
\pf
We only need to consider $(\Leftarrow)$ with the length of $P$ at least $1$. The case that $\ell \leqslant 3$ follows from Observation \ref{obs_dlr4}. Now let $\ell \geqslant 4$. For a contradiction, let $P$ be a minimal counterexample such that its ends $u, v$ are contained in two  $\ell$-paths $Q_u$ and $Q_v$ respectively. Clearly, $Q_u$ contains a sub path $L_u$ starting from $u$ and of length $s_u \geqslant \lceil\ell/2\rceil$. By the minimality of $P$, none inner vertex of $P$ belongs to $Q_u$ or $Q_v$. So the union of $L_u, P$ and $L_v$ forms a path of length at least $s_u + s_v + 1 > \ell$ in $T$, contradicting that $P$ is not $\ell$-incident.
\qed

\subsection{Incidence subgraphs}
The \emph{$\ell$-incident subgraph $G[\ell]$} of a graph $G$ is the graph induced by the $\ell$-incident units of $G$. By definition, $G[\ell] = G$ if $\ell = 0$ or $G$ is null. And for each $\ell \in [3]$, $G[\ell]$ can be obtained from $G$ by deleting all acyclic components of diameter $\leqslant \ell - 1$. For each $G \in \mathbb{R}_\ell[K_1]$, $G[\ell]$ is an $\ell$-path. The statements below, follow from the definitions and Lemma \ref{lem_cycleIncident}, allow us to concentrate on incidence subgraphs of trees of finite diameter.
\begin{coro}\label{coro_roughGl}
Let $G$ be a graph and $s \geqslant \ell \geqslant 0$. Then every $s$-link of $G$ belongs to $G[\ell]$. If further $G$ is nonnull and connected, then $G$ contains a cycle or a ray if and only if for every $\ell \geqslant 0$, $G = G[\ell]$.
\end{coro}

%

A rough structure of $T[\ell]$ can be derived from from Lemma \ref{lem_elcovered}.
\begin{coro}\label{coro_Tltree}
Let $T$ be a tree of finite $\diam(T) \geqslant \ell \geqslant 0$. Then $T[\ell]$ is an induced subtree of $T$. And each leaf of $T[\ell]$ is a leaf of $T$.
\end{coro}
\pf
By Lemma \ref{lem_elcovered}, $T[\ell] \leqslant T$. Let $v$ be a leaf of $T[\ell]$. By Corollary \ref{coro_roughGl}, there is an $\ell$-path $L$ of $T[\ell]$ with an end $v$. Suppose for a contradiction that $v$ is not a leaf of $T$. Then there exists $e \in E(T) \setminus E(T[\ell])$ incident to $v$. But $e$ and $L$ form an $(\ell + 1)$-path of $T$, contradicting that $e$ is not $\ell$-incident.
\qed

Let $u \in V(T)$ and $X$ be a subtree of $T$. Denote by $T_{X}^u$ the component of $T - E(X)$ containing $u$. Below is an accurate structure of $T[\ell]$.
\begin{lem}\label{lem_TlTu}
Let $T$ be a tree of finite $\diam(T) \geqslant \ell \geqslant 0$, and $X$ be a subtree of $T$. Then $X = T[\ell]$ if and only if $X = X[\ell]$, and for each $u \in V(X)$ either
\begin{itemize}
\item[{\bf (1)}] $\ecc_{X}(u) \geqslant \ell - 1$, and $T^u := T_{X}^u$ is a single vertex $u$. Or
\item[{\bf (2)}] $\lceil\ell/2\rceil \leqslant \ecc_{X}(u) \leqslant \ell - 2$, and $\ecc_{X}(u) + \ecc_{T^u}(u) \leqslant \ell - 1$.
\end{itemize}
\end{lem}
\pf The case of $\ell \leqslant 3$ follows from Observation \ref{obs_dlr4}. Now let $\ell \geqslant 4$. $(\Rightarrow)$ By Corollary \ref{coro_roughGl}, $s := \diam(X) = \diam(T) \geqslant \ell$, and every $u \in V(X)$ is $\ell$-incident in $X$. So $X[\ell] = X$ is nonnull, and $s \geqslant \ecc_{X}(u) \geqslant \radi(X) = \lceil s/2 \rceil \geqslant \lceil\ell/2 \rceil$. By Corollary \ref{coro_Tltree}, $T^u$ is a maximal subtree of $T$, of which the only unit that is $\ell$-incident in $T$ is the vertex $u$. Let $v \in V(T^u)$ such that $t := \dist(u, v) = \ecc_{T^u}(u)$. Then $T^u$ is not a single vertex if and only if $t \geqslant 1$. If this is the case, then $\ell - 1 \geqslant \ecc_T(v) \geqslant t + \ecc_{X}(u)$, and the statement follows. $(\Leftarrow)$ $X = X[\ell] \subseteq T[\ell]$. We still need to show that $T[\ell] \subseteq X$. Otherwise, there exists some $\vec{P} := (v_0, e_1, \ldots, e_\ell, v_\ell) \in \vec{\mathscr{L}}_\ell(T)$, and an maximum $s \in [\ell]$ such that $\vec{P}[0, s]$ belongs to $T^{v_s}$. By {\bf (2)} $\radi(T^{v_s}) \leqslant \ecc_{T^{v_s}}({v_s}) \leqslant \ell - 1 - \lceil\ell/2\rceil = \lfloor\ell/2\rfloor - 1$. So $\diam(T^{v_s}) \leqslant \ell - 2$, and hence there is a maximum $t \geqslant s + 1$ such that $t \leqslant \ell$ and $\vec{R} := \vec{P}[s, t]$ belongs to $X$. Since $\ecc_{X}(v_s) + \ecc_{T^{v_s}}(v_s) \leqslant \ell - 1$, $t \leqslant \ell - 1$. Since $X = X[\ell]$ is nonnull, there exists $\vec{L} := (\vec{L}_1.\vec{P}(s_1, t_1).\vec{L}_2) \in \vec{\mathscr{L}}_\ell(X)$, where $s \leqslant s_1 < t_1 \leqslant t$, and $L_1$ and $L_2$ are edge disjoint with $P$. Since $\ecc_{X}(v_s) + \ecc_{T^{v_s}}(v_s) < \ell$, $(\vec{P}(0, t_1).\vec{L}_2)$ is a dipath of length less than $\ell$. So $\vec{L}_2$ is of length less than $\ell - t_1$. Since $\ecc_{X}(v_t) + \ecc_{T^{v_t}}(v_t) < \ell$, $(\vec{L}_1.\vec{P}(s_1, \ell))$ is a dipath of length less than $\ell$. So $\vec{L}_1$ is of length less than $s_1$. So $\vec{L}$ is of length less than $s_1 + t_1 - s_1 + \ell - t_1 = \ell$, a contradiction.
\qed

\subsection{Equivalence classes}
In this subsection we build the relationships among $\ell$-minimal graphs, $\ell$-incident graphs, and $\ell$-equivalence classes.
\begin{lem}\label{lem_minEq}
Each $\ell$-finite graph $G$ is $\ell$-equivalent to $G[\ell]$.
\end{lem}
\pf
By Corollary \ref{coro_roughGl}, all $\ell$-links and $(\ell + 1)$-links of $G$ belong to $G[\ell] \subseteq G$. So $\mathbb{L}_\ell(G) = \mathbb{L}_\ell(G[\ell])$, and the lemma follows.
\qed

The following lemma links $\ell$-incidence units with $\ell$-minimal graphs.
\begin{lem}\label{lem_minInci}
An $\ell$-finite graph $G$ is $\ell$-minimal if and only if $G = G[\ell]$.
\end{lem}
\pf
Since every unit of $G[\ell]$ is $\ell$-incident, deleting a unit from $G[\ell]$ will delete at least one $\ell$-link from $G[\ell]$. So $G[\ell]$ is $\ell$-minimal. Conversely, if $G[\ell] \subset G$, then $G$ is not $\ell$-minimal since, by Lemma \ref{lem_minEq}, $\mathbb{L}_\ell(G[\ell]) = \mathbb{L}_\ell(G)$.
\qed

Below we connect $\ell$-equivalence relation and $\ell$-incidence graphs.
\begin{lem}\label{lem_minUniq}
For $\ell$-finite graphs $X, Y$, $X \sim_\ell Y$ if and only if $X[\ell] \cong Y[\ell]$.
\end{lem}
\pf
$(\Leftarrow)$ Let $Z := X[\ell] \subseteq X, Y$. By Lemma \ref{lem_minEq}, $\mathbb{L}_\ell(X) \cong \mathbb{L}_\ell(Z) \cong \mathbb{L}_\ell(Y)$. So $X \sim_\ell Y$. $(\Rightarrow)$ By definition there exists an $\ell$-minimal graph $Z \subseteq X, Y$ such that $\mathbb{L}_\ell(X) \cong \mathbb{L}_\ell(Y) \cong \mathbb{L}_\ell(Z)$. By Lemma \ref{lem_minInci}, $Z = Z[\ell] \subseteq X[\ell]$ since $Z \subseteq X$. But by Lemma \ref{lem_minEq}, $\mathbb{L}_\ell(Z) = \mathbb{L}_\ell(X) = \mathbb{L}_\ell(X[\ell])$. So $Z \cong X[\ell]$ since, by Lemma \ref{lem_minInci}, $X[\ell]$ is $\ell$-minimal. Similarly, $Y[\ell] \cong Z$ and the lemma follows.
\qed

\noindent {\bf Proof of Lemma \ref{lem_miniUniq}.} The reflexivity and symmetry of $\sim_\ell$ follow from the definition. To show the transitivity, let $X \sim_\ell Y$ and $Y \sim_\ell Z$. Then by Lemma \ref{lem_minUniq}, $X[\ell] \cong Y[\ell] \cong Z[\ell]$, and hence $X \sim_\ell Z$. The uniqueness of the $\ell$-minimal graph in its class follows from Lemma \ref{lem_minInci} and \ref{lem_minUniq}. The fact that $G[\ell]$ is an induced subgraph of $G$ follows from Corollary \ref{coro_roughGl} and \ref{coro_Tltree}. \qed

\noindent {\bf Proof of Lemma \ref{lem_conClass}.} Let $Z \sim_\ell G$, and $Y$ be a component of $Z$. None of $Z, Y$ and $G$ contains rays since they are $\ell$-finite. If $Y$ contains a cycle, then $Y = Y[\ell]$ is a component of $G$. Now let $Y$ be a tree. If $\diam(Y) < \ell$, then $Y[\ell] \cong K_0$. There can be arbitrarily many such $Y$. If $\radi(Y) \geqslant \ell$, $Y = Y[\ell]$ is a component of $G$.  The rest of the lemma follows from Lemma \ref{lem_TlTu}.\qed

%% file: 45partitionedGraphs.tex
\section{Partitioned $\ell$-link graphs}
A sufficient and necessary condition for an $\ell$-link graph to be connected was given by Jia and Wood \cite{JiaWood2013}. In this section, we study the cyclic components of $\ell$-roots. The investigation helps to bound the parameters of minimal $\ell$-roots.


\subsection{Definitions and basis} Let $H$ be a graph admitting partitions $\mathcal{V}$ of $V(G)$ and $\mathcal{E}$ of $E(G)$. $\tilde{H} := (H, \mathcal{V}, \mathcal{E})$ is called a \emph{partitioned graph}. For each graph $G$, let \emph{$\mathcal{V}_0(G) := \{\{v\} \subseteq V(G)\}$}, and \emph{$\mathcal{E}_0(G) := \{\{e\} \subseteq E(G)\}$}. Let $\ell \geqslant 1$. For $R \in \mathscr{L}_{\ell - 1}(G)$, let \emph{$\mathscr{L}_{\ell + 1}(R)$} be the set of $\ell$-links of $G$ of middle subsequence $R$, $\mathscr{L}^{[\ell]}(R) := \{Q^{[\ell]}| Q \in \mathscr{L}_{\ell + 1}(R)\}$, and \emph{$\mathcal{E}_\ell(G) := \{\mathscr{L}^{[\ell]}(R) \neq \emptyset | R \in \mathscr{L}_{\ell - 1}(G)\}$}. Let \emph{$\mathcal{V}_1(G) := \{M(u, v) \neq \emptyset | u, v \in V(G)\}$}. For $\ell \geqslant 2$, let \emph{$\mathcal{V}_\ell(G) := \{\mathscr{L}_\ell(R) \neq \emptyset| R \in \mathscr{L}_{\ell - 2}(G)\}$}. By \cite[Lemma 4.1]{JiaWood2013}, for $\ell \neq 1$, $\mathcal{V}_\ell(G)$ consists of independent sets of $\mathbb{L}_\ell(G)$. For $\ell \geqslant 0$, $\widetilde{\mathbb{L}}_\ell(G) := (\mathbb{L}_\ell(G), \mathcal{V}_\ell(G), \mathcal{E}_\ell(G))$ is a partitioned graph, and called a \emph{partitioned $\ell$-link graph} of $G$. $(\mathcal{V}_\ell(G), \mathcal{E}_\ell(G))$ is called an \emph{$\ell$-link partition} of $H$. $G$ is an \emph{$\ell$-root} of $\tilde{H}$ if $\widetilde{\mathbb{L}}_\ell(G) \cong \tilde{H}$. Denote by \emph{$\mathbb{R}_\ell[\tilde{H}]$} (respectively, \emph{$\mathbb{R}_\ell(\tilde{H})$}) the set of all (respectively, minimal) $\ell$-roots of $\tilde{H}$.



\begin{propo}
A graph is an $\ell$-link graph if and only if it admits an $\ell$-link partition. Moreover, $\mathbb{R}_\ell(H)$ (respectively, $\mathbb{R}_\ell[H]$) is the union of $\mathbb{R}_\ell(\tilde{H})$ (respectively, $\mathbb{R}_\ell[\tilde{H}]$) over all partitioned graph $\tilde{H}$ of $H$.
\end{propo}

An \emph{$\ell$-link} (respectively, \emph{$\ell$-arc}) of $\tilde{H}$ is an $\ell$-link (respectively, $\ell$-arc) of $H$ whose consecutive edges are in different edge parts of $\tilde{H}$. The lemma below indicates that every $s$-link of $\widetilde{\mathbb{L}}_\ell(G)$ arises from an $(\ell + s)$-link of $G$.

\begin{lem}\label{lem_sHslG}
Let $\ell, s \geqslant 0$ be integers, $G$ be a graph, and $L$ be an $s$-link of $\mathbb{L}_\ell(G)$. Then $L$ is an $s$-link of $\widetilde{\mathbb{L}}_\ell(G)$ if and only if there exists an $(\ell + s)$-link $R$ of $G$ such that $L = R^{[\ell]}$.
\end{lem}
\pf
It is trivial for $\ell = 0$ or $s \leqslant 1$. Now let $\ell \geqslant 1$, $s \geqslant 2$, and $L := [L_0, Q_1, \ldots, Q_s, L_s]$. $(\Leftarrow)$ Let $L = R^{[\ell]}$, and $P_i$ be the middle sub $(\ell - 1)$-link of $Q_i$, where $i \in [s]$. By definition $Q_i^{[\ell]} = [L_{i - 1}, Q_i, L_i] \in \mathscr{L}^{[\ell]}(P_i) \in \mathcal {E}_\ell(G)$. By \cite[Observation 3.3]{JiaWood2013}, $P_i \neq P_{i + 1}$ for $i \in [s - 1]$. So $Q_{i}^{[\ell]}$ and $Q_{i + 1}^{[\ell]}$ are in different edge parts of $\mathcal {E}_\ell(G)$. Thus $L$ is an $s$-link of $\widetilde{\mathbb{L}}_\ell(G)$. $(\Rightarrow)$ Let $L$ be an $s$-link of $\widetilde{\mathbb{L}}_\ell(G)$, where $L_i := [v_i, e_{i + 1}, \ldots, v_{i + \ell}]$ for $i \in \{0, 1\}$, and $Q_1 := [v_0, e_1, \ldots, e_{\ell + 1}, v_{\ell + 1}]$. Suppose for a contradiction that $Q_2$ has the form $[u_0, f_1, v_1, \ldots, e_{\ell + 1}, v_{\ell + 1}]$. Then $Q_2$ and $Q_1$ have the same middle sub $(\ell - 1)$-link $P := [v_1, e_2, \ldots, e_\ell, v_\ell]$, and hence are in the same part $\mathscr{L}^{[\ell]}(P) \in \mathcal {E}_\ell(G)$, a contradiction. So $Q_2$ has the form $[v_1, e_2, \ldots, e_{\ell + 2}, v_{\ell + 2}]$, and $L_2 = [v_2, e_3, \ldots, e_{\ell + 2}, v_{\ell + 2}]$. Continue this analysis, the statement follows.
\qed

As shown below, every closed $s$-link of $\widetilde{\mathbb{L}}_\ell(G)$ stems from a closed $(\ell + s)$-link of $G$.



\begin{lem}\label{lem_cycleGH}
Let $\ell \geqslant 0$, $s \geqslant 2$, $G$ be graph, and $\vec{R}$ be an $(\ell + s)$-arc of $G$. Then $R^{[\ell]}$ is a closed $s$-link of $\widetilde{\mathbb{L}}_\ell(G)$ if and only if $\vec{R}(0, \ell) = \vec{R}(s, \ell + s)$.
\end{lem}
\pf By Lemma \ref{lem_sHslG}, $[L_0, Q_1, \ldots, Q_s, L_s] := R^{[\ell]}$ is an $s$-link of $\widetilde{\mathbb{L}}_\ell(G)$. Clearly, $R^{[\ell]}$ is closed if and only if $\vec{R}[0, \ell] = L_0 = L_s = \vec{R}[s, \ell + s]$. Let $\vec{R} := (v_0, e_1, \ldots, v_{\ell + s})$. Suppose for a contradiction that $\vec{R}(0, \ell) = \vec{R}(\ell + s, s)$; that is, $(v_0, e_1, \ldots, e_\ell, v_\ell)\linebreak = (v_{\ell + s}, e_{\ell + s}, \ldots, e_{s + 1}, v_s)$. Then $Q_1 = [v_0, e_1, \ldots, e_{\ell + 1}, v_{\ell + 1}]$ and $Q_s = [v_{s - 1}, e_{s},\linebreak \ldots, e_{\ell + s}, v_{\ell + s}]$ have the same middle sub $(\ell - 1)$-link $P := [v_1, e_{2}, \ldots, e_\ell, v_\ell] = [v_{\ell + s - 1}, e_{\ell + s - 1}$, $\ldots, e_{s + 1}, v_s]$. So $Q_1, Q_s$ correspond to edges of $\widetilde{\mathbb{L}}_\ell(G)$ from the same edge part $\mathscr{L}^{[\ell]}(P)$, a contradiction. Therefore, $\vec{R}(0, \ell) = \vec{R}(s, \ell + s)$.
\qed


\subsection{Cycles in partitioned graphs}
A \emph{cycle} of $\tilde{H}$ is a cycle of $H$ whose the consecutive edges are in different edge parts. $\tilde{H}$ and its partition are \emph{cyclic} if $\tilde{H}$ contains a cycle, and \emph{acyclic} otherwise. For example, for each $t$-cycle $O$, $\widetilde{\mathbb{L}}_\ell(O)$ is cyclic. When $t \geqslant 3$ is divisible by $3$ or $4$, by Lemma \ref{lem_cyclerootT}, there exists a tree $T$ and an integer $\ell$ such that $O \cong \mathbb{L}_\ell(T)$ can be organised into an acyclic partitioned graph $\widetilde{\mathbb{L}}_\ell(T)$. Each component $X$ of $H$ corresponds to a partitioned subgraph $\tilde{X}$ of $\tilde{H}$. $\tilde{X}$ is called a \emph{component} of $\tilde{H}$. Let \emph{$o(\tilde{H})$} and \emph{$a(\tilde{H})$} be the cardinalities of the sets of cyclic and acyclic components of $\tilde{H}$ respectively. The following lemma tells that the number of cyclic components is invariant (see, for example, \cite{Olver1999}) under the partitioned $\ell$-link graph construction.

\begin{lem}\label{lem_oHeqoG}
For every graph $G$ and each integer $\ell \geqslant 0$, we have $o(G) = o(\widetilde{\mathbb{L}}_\ell(G))$.
\end{lem}
\pf
Let $\tilde{H} := \widetilde{\mathbb{L}}_\ell(G)$, $X$ be a component of $G$ containing a cycle $O$, $\tilde{X}_\ell$ be the component of $\tilde{H}$ containing $\widetilde{\mathbb{L}}_\ell(O)$. We only need to show that $\varphi: X \mapsto \tilde{X}_\ell$ is a bijection from the cyclic components of $G$ to that of $\tilde{H}$. For $i \in \{1, 2\}$, let $O_i$ be a closed $s_i$-link of $\tilde{H}$. By Lemma \ref{lem_cycleGH}, $G$ contains an $(\ell + s_i)$-link $R_i$ such that $O_i = R_i^{[\ell]}$. On one hand, if $R_i$ are closed $*$-links of $X$, by \cite[Lemma 3.3]{JiaWood2013}, $\vec{R}_i[0, \ell]$ can be shunted to each other in $X$, and the images of the shunting form a $*$-link from $O_1$ to $O_2$ in $X_\ell$. On other hand, if $O_i$ are connected by a $*$-link $Q$ of $X_\ell$ between, say, $\vec{R}_i[0, \ell] \in V(O_i)$. Then $\vec{R}_i[0, \ell] \in \mathscr{L}_\ell(G)$ can be shunted to each other, with images corresponding to the vertices of $Q$ in $X_\ell$. \qed

Every graph is a disjoint union of its components. But the relationship between a partitioned graph $\tilde{H}$ and its components is more complicated. For different components $\tilde{X}$ and $\tilde{Y}$ of $\tilde{H}$, it is possible that a vertex part $U$ of $\tilde{X}$ and a vertex part $V$ of $\tilde{Y}$ are two disjoint subsets of a vertex part $W$ of $\tilde{H}$. Lemma \ref{lem_oHeqoG} leads to a rough process of building $\widetilde{\mathbb{L}}_\ell(G)$ from its components. We can first fix all cyclic components such that no two units from different components are in the same part of $\tilde{H}$. And then for each acyclic component, we either set it as an independent fixed component, or merge some of its vertex parts with that of a unique fixed graph to get a larger fixed graph.



By definition, $o(\tilde{H}) \leqslant o(H)$, $a(H) \leqslant a(\tilde{H})$, and $a(H) + o(H) = c(H) = c(\tilde{H}) = a(\tilde{H}) + o(\tilde{H})$. As a consequence of Lemma \ref{lem_oHeqoG}, we have:
\begin{coro}\label{coro_aHeqcH} For each $\ell \geqslant 0$,
a graph $G$ is acyclic if and only if
$a(\widetilde{\mathbb{L}}_\ell(G)) = c(\mathbb{L}_\ell(G))$.
\end{coro}

A cycle of $\tilde{H}$ is \emph{plain} if it is of length $2$, or every pair of its units are in different parts. Clearly, for each cycle $O$ of $G$, $\mathbb{L}_\ell(O)$ is a plain cycle of $\widetilde{\mathbb{L}}_\ell(G)$.
\begin{obs}\label{obs_partCycle}
For each graph $G$ and every integer $\ell \geqslant 0$, a component of $\widetilde{\mathbb{L}}_\ell(G)$ contains a cycle if and only if it contains a plain cycle.
\end{obs}
\pf
By the proof of Lemma \ref{lem_oHeqoG}, a component $\tilde{X}_\ell$ of $\widetilde{\mathbb{L}}_\ell(G)$ is cyclic if and only if there is a cycle $O$ of $G$ such that $\tilde{X}_\ell$ contains the plain cycle $\mathbb{L}_\ell(O)$.
\qed

Observation \ref{obs_partCycle} helps to decide whether a partition is an $\ell$-link partition.
\begin{exa}
Let $H := [v_0, e_1, \ldots, e_6, v_6 = v_0]$ be a cycle, $\mathcal{E} := \{\{e_i\}| i \in [6]\}$, and $\mathcal {V} := \{\{v_i\},\{v_2, v_5\} | i \in \{0, 1, 3, 4\}\}$. Then $\tilde{H} := (H, \mathcal{V}, \mathcal{E})$ is cyclic with a unique closed $6$-link $H$. And $H$ is an $\ell$-link graph for every $\ell \geqslant 0$. However, $\tilde{H}$ is not a partitioned $\ell$-link graph since it contains no plain cycle.
\end{exa}

One should be aware that not every plain cycle of $\widetilde{\mathbb{L}}_\ell(G)$ is an $\ell$-link graph of a cycle of $G$. It may be formed by the sub $\ell$-links of a closed $*$-link of $G$.
\begin{exa}
$\vec{L} := (v_0, e_1, \ldots, e_3, v_0 = u_0, f_1, \ldots, f_3, u_0)$ is a closed $6$-arc but not a dicycle of a graph $G$, while $[\vec{L}.\vec{L}(0, 3)]^{[3]}$ is a plain $6$-cycle of $\widetilde{\mathbb{L}}_3(G)$.
\end{exa}





For cycles of certain lengths, we have the following:
\begin{lem}
Let $\ell \geqslant 0$, $s \geqslant 2$, $G$ be a graph, and $X$ be an $s$-cycle of $\widetilde{\mathbb{L}}_\ell(G)$. Then $g := \girth(G) \leqslant s$. If further $s \leqslant 2g - 1$. Then $G$ contains an $s$-cycle $O$ such that $X = \mathbb{L}_\ell(O)$.
\end{lem}
\pf
By Lemma \ref{lem_sHslG} and \ref{lem_cycleGH}, there exists $\vec{R} \in \vec{\mathscr{L}}_{\ell + s}(G)$ such that $X = R^{[\ell]}$, and $\vec{R}(0, \ell) = \vec{R}(s, \ell + s)$. Then $O := \vec{R}[0, s]$ is a closed $s$-link of $G$, and $s \geqslant g$. Suppose for a contradiction that $O$ is not a cycle. Without loss of generality, let $i \in [s - 1]$ such that $v_0 = v_i = v_s$; that is, $\vec{R}[0, i]$ and $\vec{R}[i, s]$ are closed $*$-links of $G$. Then $i \geqslant g$ and $s - i \geqslant g$, and hence $s \geqslant 2g$, a contradiction.
\qed

\subsection{Computing cyclic components}
We explain how to decide whether a component of $\tilde{H} := (H, \mathcal{V}, \mathcal{E})$ is cyclic, and compute $o(\tilde{H})$ efficiently.
Let \emph{$\mathcal{E}(u)$} be the subset of $\mathcal{E}$ incident to $u \in V(H)$, and \emph{$r({\mathcal{E}}) := \max\{|\mathcal{E}(u)|\; u \in V(H)\}$}.



\begin{defi}\label{defi_link2arcgraph}
Let $\tilde{H} := (H, \mathcal{V}, \mathcal{E})$ be a partitioned graph, and $\vec{H}$ be the digraph with vertices $(u, E)$, where $E \in \mathcal{E}(u)$, such that there is an arc from $(u, E)$ to $(v, F)$ if $E \neq F$, and there is $e \in E$ between $u$ and $v$.
\end{defi}




In the following, we transfer the problem of computing $o(\tilde{H})$ to that of detecting if a component of $\vec{H}$ contains a dicycle.

\begin{lem}\label{lem_p2acyclic}
In Definition \ref{defi_link2arcgraph}, let $\tilde{X}$ be a component of $\tilde{H}$. Then $\vec{X}$ consists of components of $\vec{H}$. Moreover, $\tilde{X}$ is cyclic if and only if $\vec{X}$ contains a dicycle.
\end{lem}
\pf
Let $s \geqslant 2$, $O := [v_0, e_1, \ldots, e_s, v_s := v_0]$ be a cycle of $X$, $e_i \in E_i \in \mathcal{E}$ for $i \in [s]$, $e_{s + 1} := e_1$ and $E_{s + 1} := E_1$. By definition $O$ is a cycle of $\tilde{X}$ if and only if $E_i \neq E_{i +1}$ for $i \in [s]$, which is equivalent to saying that $u_{i - 1} := (v_{i - 1}, E_i)$ is a vertex while $\vec{f}_i := (u_{i - 1}, u_i)$ is an arc of $\vec{H}$ for $i \in [s + 1]$; that is, $(u_0, \vec{f}_1, \ldots, \vec{f}_s, u_s = u_0)$ is a dicycle of $\vec{X}$.
\qed
\noindent {\bf Remark.} Let $n := n(H)$, $m := m(H)$ and $r := r(\mathcal{E})$. Then $n(\vec{H}) = \sum_{u \in V(H)}|\mathcal{E}(u)| \leqslant \sum_{u \in V(H)}\deg_H(u) = 2m$. Every arc $(u, e, v)$ of $H$ with $e \in E \in \mathcal{E}$ corresponds to $|\mathcal{E}(v)| - 1$ arcs of $\vec{H}$; that is, $((u, E), (v, F))$ for $F \in \mathcal{E}(v)\setminus \{E\}$. So $m(\vec{H}) \leqslant \sum_{(u, e, v)\in \vec{\mathscr{L}}_1(H)}(|\mathcal{E}(v)| - 1) \leqslant 2\sum_{v \in V(H)}\deg_H(v)(|\mathcal{E}(v)| - 1) \leqslant 4m(r - 1)$. An $O(n + m)$-time algorithm for dividing $H$ into connected components was given by Hopcroft and Tarjan \cite{HopcroftTarjan1973}. For each component $\tilde{X}$ of $\tilde{H}$, Tarjan's algorithm \cite{Tarjan1972}, with time complexity $O(n(\vec{X}) + m(\vec{X})) = O(rm(X))$, can be used to detect the existence of dicycles in $\vec{X}$, and hence that of cycles in $\tilde{X}$ by Lemma \ref{lem_p2acyclic}. So the time complexity for computing $o(\tilde{H})$ is $O(m^2 + n)$ in general situations, and is $O(m + n)$ if $r$ is bounded. By \cite[Lemma 4.1]{JiaWood2013}, $r(\mathcal{E}_\ell(G)) \leqslant 2$ for $\ell \geqslant 1$. So it requires $O(m + n)$-time to compute the number of cyclic components of a partitioned $\ell$-link graph for each $\ell \geqslant 0$.

%% file: 46BoundingSize.tex
\section{Bounding the number of minimal roots}\label{sec_boundSize}
We bound in this section the order, size, maximum degree and total number of minimal $\ell$-roots of a finite graph. This lies the basis for solving the recognition and determination problems for $\ell$-link graphs in our future work.



\subsection{Incidence pairs}
Let $s, \ell \geqslant 0$, $L \in \mathscr{L}_\ell(G)$, \emph{$\mathcal {I}_G(L, s)$} be the set of pairs $(L, R)$ such that $R \in \mathscr{L}_s(G)$ is incident to $L$, and \emph{$i_G(L, s) := |\mathcal {I}_G(L, s)|$}. \emph{$\girth(L)$} is $+\infty$ if $L$ is a path, and the minimum length of a sub cycle of $L$ otherwise. To dodge confusions, denote by \emph{$\hat{L}$} the graph induced by the units of $L$. Then $\girth(\hat{L}) \leqslant \girth(L)$ and the inequation may hold. For example, when $L= [v_0, e_1, v_1, e_2, v_2, e_3, v_0, e_4, v_1]$, $\girth(\hat{L}) = 2$ while $\girth(L) = 3$.
\begin{lem}\label{lem_inciCounting}
Let $\ell \geqslant s \geqslant 0$, $L \in \mathscr{L}_\ell(G)$, and $g := \girth(L^{[s]})$. Then $\min\{g, \ell - s + 1\} \leqslant i_G(L, s) \leqslant \ell - s + 1$. Further, $i_G(L, s) = \ell - s + 1$ if and only if $g = +\infty$ if and only if $L^{[s]}$ is an $(\ell - s)$-path if and only if $g \geqslant \ell - s + 1$. Otherwise, if $g \leqslant \ell - s$, then $i_G(L, s) = g$ if and only if $\hat{L^{[s]}}$ is a $g$-cycle.
\end{lem}
\pf
$L^{[s]}$ is an $(\ell - s)$-link of $H := \mathbb{L}_s(G)$. So $t := i_G(L, s) =  |L^{\{s\}}| \leqslant \ell - s + 1$, with the last equation holds if and only if $L^{[s]}$ is a path. If $g \leqslant \ell - s$, then $L^{[s]}$ contains a sub $g$-cycle on which every pair of different vertices of $H$ corresponding to a pair of different sub $s$-links of $L$. So $t \geqslant g$, with equation holds if and only if all units of $H$ on $L^{[s]}$ belong to the $g$-cycle.
\qed


Let \emph{$\mathcal {I}_G(\ell, s) := \bigcup_{L \in \mathscr{L}_\ell(G)}\mathcal {I}_G(L, s)$}, and \emph{$i_G(\ell, s):= |\mathcal {I}_G(\ell, s)|$}. Then $i_G(\ell, s) = \sum_{L \in \mathscr {L}_\ell(G)}$ $i_G(L, s)$ can be bounded as follows.


\begin{coro}\label{coro_inciCounting}
Let $\ell \geqslant s \geqslant 0$, and $G$ be $\ell$-finite of girth $g$. Then
\begin{eqnarray*}\label{eq_gllink}
\min\{g, \ell - s + 1\}|\mathscr{L}_\ell(G)| &\leqslant& i_G(\ell, s) \leqslant (\ell - s + 1)|\mathscr{L}_\ell(G)|.
\end{eqnarray*}
Further, $i_G(\ell, s) = (\ell - s + 1)|\mathscr{L}_\ell(G)|$ if and only if $g \geqslant \ell - s + 1$. If $g \leqslant \ell - s$, then $i_G(\ell, s) = g|\mathscr{L}_\ell(G)|$ if and only if $G$ is a disjoint union of $g$-cycles.
\end{coro}
\pf
Let $t: = \min\{\girth(L^{[s]})| L\in \mathscr{L}_\ell(G)\}$. Note that $g \geqslant \ell - s + 1$ if and only if $t = +\infty$. And $g \leqslant \ell - s$ if and only if $t = g$. So the statements can be verified by summing the results in Lemma \ref{lem_inciCounting} over $L \in \mathscr{L}_\ell(G)$.
\qed

The order and size of minimal $\ell$-roots are bounded as follows.
\begin{lem}\label{lem_HboundG}
Let $\ell \geqslant 1$ and $G$ be an $\ell$-minimal graph. Then both $G$ and $\tilde{H} := \widetilde{\mathbb{L}}_\ell(G)$ are finite. Further, $m(G) \leqslant \ell n(H)$, and $n(G) \leqslant \ell n(H) + a(\tilde{H})$.
\end{lem}
\pf
Since $G$ is $\ell$-finite, $\tilde{H}$ is finite. By Corollary \ref{coro_inciCounting}, $i_G(\ell, 1) \leqslant \ell n(H)$ is finite. By Lemma \ref{lem_minInci}, for each $e \in E(G)$, $i_G(e, \ell) \geqslant 1$. Summing this inequation over $e \in E(G)$, we have $m(G) \leqslant i_G(\ell, 1)$ is finite. By Lemma \ref{lem_oHeqoG}, $o(G) = o(\tilde{H})$, and hence $a(G) = c(G) - o(G) \leqslant c(\tilde{H}) - o(\tilde{H}) = a(\tilde{H})$.
So $n(G) \leqslant m(G) + a(G) \leqslant \ell n(H) + a(\tilde{H})$ is finite.
\qed

\subsection{The number of minimal roots} We have known that $|\mathbb{R}_0(G)| = 1$, and $|\mathbb{R}_\ell(K_s)| = 1$ for $s \in \{0, 1, 2\}$. By Lemma \ref{lem_barK2}, $|\mathbb{R}_\ell(\bar{K}_2)| = \lfloor\frac{\ell + 1}{2}\rfloor$ for $\ell \geqslant 1$. Let $\tilde{a}_\ell(H)$ be the maximum $a(\tilde{H})$ over partitioned $\ell$-link graph $\tilde{H} := (H, \mathcal{V}, \mathcal{E})$. Clearly, $\tilde{a}_\ell(H) \leqslant c(H)$. Denote by \emph{$\psi(p, q)$} the number of nonisomorphic graphs with $p$ vertices and $q$ edges. Then $\psi(p, q)$ is $0$ if $p \leqslant 1$ and $q \geqslant 1$; is $1$ if $p = 2$ or $q = 0$; is $1$ if $p \geqslant 2$ and $q \leqslant 1$; is less than $\binom{p}{2}^q$ if $p \geqslant 3$ and $q \geqslant 1$; and by Hardy \cite{Hardy1920}, is the nearest integer to $\frac{(q + 3)^2}{12}$ if $p = 3$.
Further, when $p \geqslant 3$ and $q \geqslant 1$, $\binom{p}{2} \geqslant 3$, $\binom{p - 1}{2}/\binom{p}{2} = (1 - \frac{2}{p})^{\frac{p}{2}\cdot\frac{2}{p}} < \textbf{e}^{\frac{-2}{p}}$, and so $\sum_{i = 0}^p\sum_{j = 0}^q \psi(i, j) \leqslant p + q + \sum_{i = 3}^p\sum_{j = 1}^q\binom{i}{2}^j < p + q + \frac{3}{2}\sum_{i = 3}^p\binom{i}{2}^q < p + q + \frac{3}{2}\binom{p}{2}^q[1 - \textbf{e}^{\frac{-2q}{p}}]^{-1}$.

\begin{lem}\label{lem_lminiFinite}
Let $\ell \geqslant 1$ be an integer, $H$ be a finite graph, $n := n(H) \geqslant 2$, and $a := \tilde{a}_\ell(H)$. Then $H$ has at most $(\ell n + a)^{2\ell n}$ minimal $\ell$-roots, of which at most $2(\ell n + 1)^{\ell n - 1}$ are trees, and at most $a(a + 1)(\ell n + a)^{\ell n - 1}$ are forests.
\end{lem}
\pf Let $G \in \mathbb{R}_\ell(H)$. By Lemma \ref{lem_HboundG}, $m(G) \leqslant \ell n$, and $n(G) \leqslant \ell n + a \leqslant 2\ell n$. By the analysis above, $|\mathbb{R}_\ell(H)| < 2\ell n + a + \frac{3}{2}\binom{\ell n + a}{2}^{\ell n}(1 - \textbf{e}^{\frac{-2\ell n}{\ell n + a}})^{-1} \leqslant 2\ell n + a + \frac{3\textbf{e}}{2(\textbf{e} - 1)}\binom{\ell n + a}{2}^{\ell n} < (\ell n + a)^{2\ell n}$.


Cayley's formula \cite{Borchardt1860,Cayley2009} states that there are $p^{p - 2}$ unequal trees on vertex set $[p]$. So the number of trees in $\mathbb{R}_\ell(H)$ is at most $\sum_{p = 1}^{\ell n + 1}p^{p - 2} \leqslant 2(\ell n + 1)^{\ell n - 1}$.

By Aigner and Ziegler \cite{AignerZiegler2010}, the number of unequal forests on vertex set $[p]$ of $k$ components is $kp^{p - k - 1}$. So the number of forests in $\mathbb{R}_\ell(H)$ is at most $\sum_{k = 1}^a\sum_{p = a}^{\ell n + k}kn^{p - k - 1} \leqslant 2\sum_{k = 1}^{a}k(\ell n + k)^{\ell n - 1} \leqslant a(a + 1)(\ell n + a)^{\ell n - 1}$.
\qed

%% file: 47BoundingMaximumDegree.tex
\subsection{The maximum degree of minimal roots} We have proved that minimal $\ell$-roots of a finite graph are finite. So in this subsection, we only deal with finite graphs. Let $\mathcal{E}$ be a partition of $E(H)$. We may identify $E \in \mathcal{E}$ with the subgraph of $H$ induced by $E$. Let \emph{$D(\mathcal {E})$} be the set of $\deg_E(v)$ over $E \in \mathcal{E}$ incident to $v \in V(H)$. Let \emph{$D(G) := \{\deg_G(v) - 1 \geqslant 1 | v \in V(G)\}$}. Clearly, $D(G) = \emptyset$ if and only if $\Delta(G) \leqslant 1$ if and only if for all $\ell \geqslant 1$, $D(\mathcal{E}_\ell(G)) = \emptyset$.

\begin{lem}\label{lem_Dsdecrease}
Let $\ell \geqslant 1$, $G$ be a finite connected graph of $\Delta(G) \geqslant 2$, $D := D(G)$, and $D_\ell := D(\mathcal {E}_\ell)$. If $G$ is cyclic, $D = D_\ell$. Otherwise, $G$ is a tree of diameter $s \geqslant 2$, and $D = D_1 \supseteq \ldots \supseteq D_{s - 1} \supset \emptyset = D_{s} = D_{s + 1} = \ldots$.
\end{lem}
\pf For $d, \ell \geqslant 1$, by definition $d \in D_\ell$ if and only if there is an $\ell$-arc $\vec{L}$ starting from some $v \in V(G)$ with $\deg_G(v) = d + 1$. So $D = D_1$. When $\ell \geqslant 2$, $\vec{L}(0, \ell - 1)$ is an $(\ell - 1)$-arc starting from $v$. So $d \in D_{\ell - 1}$. In another word, $D_{\ell - 1} \supseteq D_\ell$ for all $\ell \geqslant 2$. On one hand, let $G$ be a tree of diameter $s \geqslant 2$. By definition $\ell \geqslant s$ if and only if $E(\mathbb{L}_\ell(G)) = \emptyset$ if and only if $D_\ell = \emptyset$. On other hand, let $O$ be a cycle of $G$. For each $v \in V(G)$, since $G$ is connected, there is a dipath $\vec{P}$ of minimum length from $v$ to some $u \in V(O)$. Clearly, there is an $\ell$-arc $\vec{R}$ of $O$ starting from $u$. Then $\vec{L} := (\vec{P}.\vec{R})(0, \ell)$ is an $\ell$-arc of $G$ starting from $v$. Thus $D_\ell = D$ by the analysis above.
\qed

Define \emph{$\Delta(\mathcal {E})$} to be $0$ if $\mathcal {E} = \emptyset$, and $\max(D(\mathcal {E}))$ otherwise. A lower bound of $\Delta(G)$ for $G \in \mathbb{R}_\ell[\tilde{H}]$ follows from Lemma \ref{lem_Dsdecrease}:


\begin{coro}\label{coro_DsdecreaseCyc}
Let $\ell \geqslant 1$, $G$ be a finite graph, and $\mathcal{E} := \mathcal{E}_\ell(G)$. Then $\Delta(G) = \Delta(\mathcal{E}) = 0$ if and only if $E(G) = \emptyset$. And $E(G) \neq \emptyset$ if and only if $\Delta(G) \geqslant \Delta(\mathcal{E}) + 1$. Further assume that $G$ is connected. Then $G$ contains a cycle if and only if for all $\ell \geqslant 1$,
$\Delta(G) = \Delta(\mathcal{E}) + 1 \leqslant \Delta(H)$.
\end{coro}

Below we display a connection between the degrees of an $\ell$-minimal tree and the number of components of the $\ell$-link graph of the tree.

\begin{lem}\label{lem_Tdc}
Let $\ell \geqslant 1$ be an integer, $T$ be a finite $\ell$-minimal tree, and $v$ be a vertex of eccentricity less than $\ell$ in $T$. Then $\deg_T(v) \leqslant c(\mathbb{L}_\ell(T)) + 1$.
\end{lem}
\pf
Clearly, $s := \ecc_T(v) \geqslant 1$ and $\ell \geqslant s + 1 \geqslant 2$. If $d := \deg_T(v) \leqslant 1$, there is nothing to show. Now let $d \geqslant 2$. For $i \in [d]$, let $\vec{e}_i := (v, e_i, u_i)$ be the arcs of $T$ starting from $v$. Then there exists $\vec{R} \in \vec{\mathscr{L}}_s(T)$ starting from, say, $\vec{e}_d$. Since $T$ is $\ell$-minimal, $e_i$ is $\ell$-incident in $T$. For $i \in [d - 1]$, by Lemma \ref{lem_elcovered}, $t_i := \ecc_{T_{e_i}^{u_i}}(u_i) \geqslant \ell - s - 1 \geqslant 0$. So there is a $t_i$-arc $\vec{Q}_i$ from $u_i$ to some $v_i$ in $T_{e_i}^{u_i}$. Obviously, $L_i := [\vec{Q}_i(\ell - s - 1, 0). \vec{e}_i. \vec{R}]$, for $i \in [d - 1]$, are $d - 1$ different $\ell$-paths containing $v$ in $T$. Suppose for a contradiction that $L_i$ can be shunted to $L_j$ for some $1 \leqslant i < j \leqslant d - 1$. Since $e_i$ separates $v$ from $v_i$, $v$ is an image of $v_i$ during the shunting. Then $\ecc_T(v) \geqslant \ecc_T(v_1) \geqslant \ell$, a contradiction. So $L_i$ and $L_j$ correspond to vertices in different components of $H := \mathbb{L}_\ell(T)$. Hence $d - 1 \leqslant c(H)$.
\qed

We now bound the maximum degree of a finite tree in terms of $\widetilde{\mathbb{L}}_\ell(T)$.

\begin{coro}\label{coro_DeltaT}
Let $\ell \geqslant 1$, $T$ be a finite tree, $\tilde{H} := (H, \mathcal {V}, \mathcal {E}) := \widetilde{\mathbb{L}}_\ell(T)$ and $s := \max\{\ecc_T(v) |$ $\deg_T(v) = \Delta(T)\}$. Then
\begin{itemize}
\item[{\bf (1)}] If $s \geqslant \ell + 1$, then $\Delta(T) = \Delta(\mathcal {E}) + 1 \leqslant \Delta(H)$.
\item[{\bf (2)}] If $\ell = s$, then either $s  = \Delta(T) = 1$, and $\Delta(H) = 0$; or $s \geqslant 2$ and $\Delta(T) = \Delta(\mathcal {E}) + 1 \leqslant \Delta(H) + 1$.
\item[{\bf (3)}] If $\ell \geqslant s + 1$ and $T$ is $\ell$-minimal, then $\Delta(T) \leqslant a(\tilde{H}) + 1 \leqslant c(H) + 1$.
\end{itemize}
\end{coro}
\pf
{\bf (1)} and {\bf (2)} are implied by the proof of Lemma \ref{lem_Dsdecrease}. {\bf (3)} follows from Lemma \ref{lem_Tdc} and Corollary \ref{coro_aHeqcH}.
\qed


For $\tilde{H} := (H, \mathcal{V}, \mathcal{E})$, let $b(\tilde{H}) := \max\{a(\tilde{H}), \Delta(\mathcal{E})\} \leqslant \max\{c(H), \Delta(H)\}$. Below we  bound the maximum degree of minimal $\ell$-roots $G$ of $\tilde{H}$ or $H$. The results in turn helps to bound the parameters of $s$-link graphs of $G$. Clearly, for $s \geqslant 1$, a graph of $q \geqslant 2$ edges contains at most $q(q - 1)^{s - 1}$ $s$-links, with equation holds if and only if all these edges are between the same pair of vertices.
\begin{lem}
Let $\ell, s \geqslant 1$, $n := n(H) \geqslant 2$, $b := b(\tilde{H})$, $G \in \mathbb{R}_\ell(\tilde{H})$, and $(X, \mathscr{V}, \mathscr{E}) := \widetilde{\mathbb{L}}_s(G)$. Then $\Delta(G) \leqslant b + 1$, $\Delta(\mathscr{E}) \leqslant b$, $\Delta(X) \leqslant 2b$, $\max\{|V|, |E| \\ |\, V \in \mathscr{V}, E \in \mathscr{E}\} \leqslant b^2$, $m(X) \leqslant \ell n(\ell n - 1)^s$, and $n(X) \leqslant \ell n(\ell n - 1)^{s - 1}$.
\end{lem}
\pf
Let $Y, Z$ be the subgraphs induced by the cyclic and acyclic components of $G$ respectively. Let $(H, \mathcal{V}, \mathcal{E}) := \tilde{H}$. By Corollary \ref{coro_DsdecreaseCyc}, $\Delta(Y) \leqslant \Delta(\mathcal{E}) + 1$. By Corollary \ref{coro_DeltaT}, $\Delta(Z) \leqslant a(\tilde{H}) + 1$. So $\Delta(G) = \max\{\Delta(Y), \Delta(Z)\} \leqslant b + 1$. The rest of the lemma follows from the analysis above.
\qed

\subsection{Better-quasi-ordering of roots}
For $i \in \{1, 2\}$, let $T_i$ be a rooted tree of root $v_i$. Denote by \emph{$T_1 \leqslant T_2$} that there is an isomorphism $\varphi$ from $T_1$ to a subtree of $T_2$ such that $\varphi(v_1) = v_2$. Clearly, $\leqslant$ defines a quasi-ordering on the set of rooted trees.
It is natural to extend $\leqslant$ to the power set of rooted trees. Let $F_1$ and $F_2$ be two sets of rooted trees. Write \emph{$F_1 \leqslant F_2$} if there is an injection $\sigma$ from $F_1$ to $F_2$ such that for every $T \in F_1$, $T \leqslant \sigma(T)$.

\noindent {\bf Proof of Lemma \ref{lem_treeBetterInsub}.}
Let $\mathcal{T}_\ell$ be the set of rooted trees of height $\ell$ ordered by $\leqslant$ defined above. By Lemma \ref{lem_coverbqo},
it is enough to show that $\mathcal{T}_\ell$ is better-quasi-ordered. The case of $\ell = 0$ follows from Lemma \ref{lem_subseqbqo}. Inductively assume it holds for some $\ell - 1 \geqslant 0$. Then by Lemma \ref{lem_subseqbqo}, the powerset $\mathcal{F}_{\ell - 1}$ of $\mathcal{T}_{\ell - 1}$ is better-quasi-ordered. For each rooted tree $T$ of root $v$, let $T(v)$ be the set of rooted trees that are components of $T - v$ with roots the neighbours of $v$ in $T$.
Note that $T \mapsto T(v)$, for $T \in \mathcal{T}_{\ell}$, is an order preserving bijection from $\mathcal{T}_{\ell}$ to $\mathcal{F}_{\ell - 1}$. So $\mathcal{T}_{\ell}$ is better-quasi-ordered since $\mathcal{F}_{\ell - 1}$ is better-quasi-ordered.\qed

\noindent {\bf Proof of Theorem \ref{thm_llinkbqo}.} Lemma \ref{lem_lminiFinite} ensures that $\mathbb{R}_\ell[H]$ is the union of finitely many $\ell$-equivalence classes. By Lemma \ref{lem_coverbqo}, it is enough to show that every $\ell$-equivalence class $\mathcal{R}$ is better-quasi-ordered by the induced subgraph relation. Let \emph{$O(G)$} (respectively, \emph{$A(G)$}) be the subgraph induced by the cyclic (respectively, acyclic) components of a graph $G$. By Lemma \ref{lem_conClass}, $\mathcal{O} := \{O(G) | G \in \mathcal{R}\}$ contains at most one graph up to isomorphism, and hence is better-quasi-ordered by Lemma \ref{lem_subseqbqo}. For each component $T$ of $A(G)$, by Lemma \ref{lem_HboundG}, $\diam(T) \leqslant m(T) \leqslant \ell n(H)$. Thus by Lemma \ref{lem_treeBetterInsub} and \ref{lem_coverbqo}, $\mathcal{A} := \{A(G) | G \in \mathcal{R}\} \subseteq \mathcal{F}_{\ell n(H)}$ is better-quasi-ordered, and so is $\mathcal{O} \times \mathcal{A}$ by Lemma \ref{lem_subseqbqo}. Note that $G \mapsto O(G) \times A(G)$, for $G \in \mathcal{R}$, is an order-preserving injection from $\mathcal{R}$ to $\mathcal{O} \times \mathcal{A}$. So by Lemma \ref{lem_coverbqo}, $\mathcal{R}$ is better-quasi-ordered.\qed

%% file: 48PathGraphs.tex
\section{Path graphs}
Some ideas and techniques used in the investigation of $\ell$-link graphs can be transplanted to study $\ell$-path graphs. We end this paper by bounding the parameters of minimal $\ell$-path roots, and showing the better-quasi-ordering results for the $\ell$-path roots of finite graphs.
\subsection{Quantitative analysis}
A graph is \emph{$\ell$-path finite} if its $\ell$-path graph is finite. Note that an $\ell$-finite graph is $\ell$-path finite, but not vise versa. For example, the disjoint union of infinitely many $\ell$-cycles is $\ell$-path finite but not $\ell$-finite.
Two $\ell$-path finite graphs $X$ and $Y$ are \emph{$\ell$-path equivalent}, written $X \simeq_\ell Y$, if there exists some graph $Z \subseteq X, Y$ such that $\mathbb{P}_\ell(X) \cong \mathbb{P}_\ell(Y) \cong \mathbb{P}_\ell(Z)$. An $\ell$-path finite graph $G$ is said to be \emph{$\ell$-path minimal} if either $G$ is null, or for each $X \subset G$, $\mathbb{P}_\ell(X) \subset \mathbb{P}_\ell(G)$. Similar with Lemma \ref{lem_miniUniq}, we have:
\begin{propo}
$\simeq_\ell$ defines an equivalence relation on $\ell$-path finite graphs. Further, each $\ell$-path equivalence class contains a unique (up to isomorphism) minimal graph, which is $\ell$-path minimal. Moreover, an $\ell$-path minimal graph is a subgraph of every graph in its $\ell$-path equivalence class.
\end{propo}

We exemplify as follows that an $\ell$-minimal graph may not be $\ell$-path minimal. And an $\ell$-path minimal graph may not be an induced subgraph of another graph in its $\ell$-path equivalence class.
\begin{exa}
Let $G$ be a graph obtained from a path $P = [v_0,\ldots, v_4]$ by adding an edge $e$ between $v_1$ and $v_3$. By Lemma \ref{lem_cycleIncident} and \ref{lem_minInci}, $G$ and $P$ are $4$-minimal. Since $\mathbb{P}_4(G) = \mathbb{P}_4(P) \cong K_1$, $G$ is not $\ell$-path minimal. Moreover, $P$ is $4$-path minimal and is a subgraph but not an induced subgraph of $G$.
\end{exa}



%
%

The following example indicates that, a graph can be both $\ell$-minimal and $\ell$-path minimal. But even so its $\ell$-link and $\ell$-path graphs may not be isomorphic.
\begin{exa}
The complete graph $K_{\ell + 1}$ is both $\ell$-minimal and $\ell$-path minimal. Clearly, for $\ell = 0, 1$ and $2$, $\mathbb{L}_\ell(K_{\ell + 1}) = \mathbb{P}_{\ell}(K_{\ell + 1})$ is isomorphic to $K_1, K_1$ and $K_3$ respectively. Now let $\ell \geqslant 3$. By \cite[Theorem 1.2]{JiaWood2013}, $\mathbb{L}_\ell(K_{\ell + 1})$ contains a $K_{\ell + 1}$-minor. However, $\mathbb{P}_{\ell}(K_{\ell + 1})$ consists of $\ell!/2$ disjoint cycles of length $\ell + 1$.
\end{exa}


%
%

By definition $\mathbb{Q}_0(H) = \{H\}$ if $H$ is simple.  $\mathbb{Q}_\ell(K_0) = \{K_0\}$, $\mathbb{Q}_\ell(K_1)$ consists of an $\ell$-path. $\mathbb{Q}_1(K_2)$ consists of a $2$-path and a $2$-cycle. When $\ell \geqslant 2$, $\mathbb{Q}_\ell(K_2)$ consists of an $(\ell + 1)$-path. Similar with Lemma \ref{lem_lminiFinite}, the order, size and the total number of minimal $\ell$-path roots of a finite graph are bounded.
\begin{lem}\label{lem_pathrootosnbound}
Let $\ell \geqslant 1$, $H$ be a finite graph, $c := c(H)$, $n := n(H) \geqslant 2$, and $G \in \mathbb{Q}_{\ell}(H)$. Then $n(G) \leqslant \ell n + c$ and $m(G) \leqslant \ell n$, with each equation holds if and only if  $G$ is a disjoint union of $\ell$-paths. Moreover, $\mathbb{Q}_{\ell}(H)$ contains at most $(\ell n + c)^{2\ell n}$ graphs, in which at most $2(\ell n + 1)^{\ell n - 1}$ are trees, and at most $c(c + 1)(\ell n + c)^{\ell n - 1}$ are forests.
\end{lem}

\subsection{$\ell$-path free graphs}
A graph is \emph{$\ell$-path free} if it contains no $\ell$-path as a subgraph, or equivalently, as a minor. By considering the \emph{types} of graphs, Ding \cite[Theorem 2.2]{Ding1992} proved that finite simple $\ell$-path free graphs are well-quasi-ordered by induced subgraph relation. Another proof, based on \emph{tree-depth} was given by Ne{\v{s}}et{\v{r}}il and Ossona de Mendez in \cite[Lemma 6.13]{NesetrilOssona2012}. We generalise this result to infinite graphs, possibly with parallel edges and loops.




\begin{lem}\label{lem_v01vl1l}
Let $\ell \geqslant 1$, $G$ be a graph with or without loops, $(T, \mathcal{V})$ be a tree-decomposition of $G$, and $[v_0, e_1, \ldots, e_\ell, v_\ell]$ be a path of $T$. Then $U = V_{v_0}\cap V_{v_1} = V_{v_{\ell - 1}}\cap V_{v_\ell}$ if and only if $U = V_{v_i} \cap V_{v_0} = V_{v_i} \cap V_{v_\ell}$ for all $i \in [\ell - 1]$.
\end{lem}
\pf $(\Leftarrow)$ is straightforward. For $(\Rightarrow)$, clearly $U = V_{v_0} \cap V_{v_1} \subseteq V_{v_0} \cap V_{v_\ell}$. By the definition of tree-decompositions, $V_{v_0} \cap V_{v_\ell} \subseteq V_{v_0} \cap V_{v_i} \subseteq V_{v_0} \cap V_{v_1}$. So $U = V_{v_0} \cap V_{v_i}$. Similarly, $U = V_{v_i} \cap V_{v_\ell}$, and hence the lemma follows.
\qed

\begin{lem}\label{lem_reducelength}
In Lemma \ref{lem_v01vl1l}, let $U = V_0 \cap V_{1} = V_{\ell - 1} \cap V_{\ell}$ for some $\ell \geqslant 3$, and $\mathcal{V}^* := \mathcal{V} \cup \{U\}$. If $U \in \mathcal{V}$, let $\mathcal{V}_u := U$ and $T' = T$. Otherwise, create a new node $u$ corresponding to $U$, and let $T'$ be obtained from $T$ by replacing some edge $e_{k + 1}$ by a $2$-path of middle vertex $u$, where $k \in [\ell - 2]$. Let $T^*$ be obtained from $T'$ by: for $i \in \{0, \ell\}$, if $u \neq v_{|i - 1|}$, then add an edge between $u$ and $v_i$, and delete the edge between $v_{|i - 1|}$ and $v_i$.
Then $(T^*, \mathcal{V}^*)$ is a tree-decomposition of $G$ of width $\tw(T, \mathcal{V})$. Further, $(T^*, \mathcal{V}^*)$ is linked if and only if $(T, \mathcal{V})$ is linked.
\end{lem}
\pf Clearly, $(T, \mathcal{V})$ is a (respectively, linked) tree-decomposition if and only if $(T', \mathcal{V}^*)$ is a (respectively, linked) tree-decomposition. And if this is the case, $\tw(T, \mathcal{V}) = \tw(T', \mathcal{V}^*)$. Obviously, $T^*$ is a tree, $\bigcup\mathcal{V}^* = V(G)$, and every pair of adjacent vertices are contained in some $V \in \mathcal{V}^*$. It follows from Lemma \ref{lem_v01vl1l} that, for different $v, w \in V(T') = V(T^*)$, ${\mathcal{V}^*}_{T'}(v, w) = {\mathcal{V}^*}_{T^*}(v, w)$. So $(T', \mathcal{V}^*)$ is a (respectively, linked) tree-decomposition if and only if $(T^*, \mathcal{V}^*)$ is a (respectively, linked) tree-decomposition. \qed

A tree-decomposition $(T, \mathcal{V})$ is \emph{short} if it has no repeated parts, and for each path $[v_0, \ldots, v_\ell]$ of $T$, where $\ell \geqslant 3$, $V_{v_0} \cap V_{v_1} \neq V_{v_{\ell - 1}} \cap V_{v_{\ell}}$. An \emph{$M$-closure} of $G$ is a triple $(T, \mathcal{V}, X)$, where $X$ is a chordal graph without a complete subgraph of order $\tw(G) + 2$, $V(G) = V(X)$, $E(G) \subseteq E(X)$, $(T, \mathcal{V})$ is a linked tree-decomposition of $X$ such that each part induces a maximal complete subgraph of $X$. An $M$-closure is \emph{short} if the tree-decomposition is short.

\begin{lem}\label{lem_linkDshort}
Every graph $G$ of finite tree-width, with or without loops, admits a short linked tree-decomposition of width $\tw(G)$.
\end{lem}
\pf By K{\v{r}}{\'{\i}}{\v{z}} and Thomas \cite[(2.3)]{KrizThomas1991}, every finite graph has an $M$-closure $(T, \mathcal{V}, X)$. Repeat the operations in Lemma \ref{lem_reducelength} on $(T, \mathcal{V})$ until no path of $T$ of length at least $3$ satisfies the conditions in Lemma \ref{lem_v01vl1l}. Delete the repeated parts, we obtain a short $M$-closure. In \cite[(2.4)]{KrizThomas1991}, replacing `an $M$-closure' by `a short $M$-closure' causes no conflict. So every graph has a short $M$-closure. The rest of the lemma follows from a discussion similar with \cite[(2.2)]{KrizThomas1991}.
\qed


\begin{obs}\label{obs_tdiDiscConn}
Let $G$ be a graph, and $\tdi(X)$ be the maximum over all components $X$ of $G$. Then $\tdi(G) \leqslant \tdi(X) + 3$.
\end{obs}

Below we show that the tree-diameter is bounded for $\ell$-path free graphs.

\begin{lem}\label{lem_tdibounded}
Every $\ell$-path free graph $G$ admits a linked tree-decomposition of width at most $\ell - 1$, and diameter  at most $2(\ell^2 - \ell + 2)^\ell + 1$. If further $G$ is connected, then $\tdi(G) \leqslant 2(\ell^2 - \ell + 2)^\ell - 2$.
\end{lem}
\pf
By Robertson and Seymour \cite{RobertsonSeymour1983}, for each finite subgraph $X$ of $G$, if $\tw(X) \geqslant \ell$, then $X$ contains an $\ell$-path, a contradiction. So $\tw(X) \leqslant \ell - 1$. By a compactness theorem for the notion of tree-width \cite{Thomas1988,Thomassen1989}, $\tw(G) \leqslant \ell - 1$.

For the tree-diameter, by Observation \ref{obs_tdiDiscConn}, we only need to consider the case that $G$ is nonnull and connected. By Lemma \ref{lem_linkDshort}, $G$ admits a short linked tree-decomposition $(T, \mathcal{V})$ of width $\tw(G)$. Let $p := \tw(G) + 1 \geqslant 1$. For $s \geqslant 1$ and $t \in [s]$, $P := [v_0, e_1, \ldots, e_s, v_s] \in \mathscr{P}_s(T)$ is called \emph{$t$-rotund} if there exists $k \in [p]$ and $1 \leqslant i_1 < \ldots < i_t \leqslant s$ such that $V_{e_{i_1}}, \ldots, V_{e_{i_t}}$ are pairwise distinct, $|V_{e_{i_j}}| = k$ for $j \in [t]$, and  $|V_{e_{j}}| \geqslant k$ for $i_1 \leqslant j \leqslant i_t$. Let $s^* \in [s]$ be the maximum number of edges of $P$ corresponding to pairwise different vertex parts of $G$.

\noindent {\bf Claim.} $s \leqslant 2s^*$. Since otherwise, there are $1 \leqslant j_1 < j_2 < j_3 \leqslant s$ such that $V_{e_{j_1}} = V_{e_{j_2}} = V_{e_{j_3}}$, contradicting the shortness of $(T, \mathcal{V})$.

\noindent {\bf Claim.} If $P$ is not $t$-rotund, then $s^* \leqslant t^p - 1$. To see this, let $s_k := |\{j |\; |V_{e_{i_j}}| = k\}|$. Since $P$ is not $t$-rotund, $s_1 \leqslant t - 1$, and for $k \geqslant 2$, $s_k \leqslant (s_1 + \ldots + s_{k - 1} + 1)(t - 1)$. An induction on $k$ gives $s_k \leqslant t^{k - 1}(t - 1)$ for each $k \in [p]$. So $s^* = s_1 + \ldots + s_p \leqslant t^p - 1$.

\noindent {\bf Claim.} If $P$ is $t$-rotund and $G$ is $\ell$-path free, then $t \leqslant p(\ell - 1) + 1$. To prove this, recall that we are considering linked tree-decomposition. So there are $k$ disjoint paths in $G$ with at least $|\bigcup_{j = 1}^t V_{e_{i_j}}| \geqslant k + t - 1$ vertices. Since $G$ is $\ell$-path free, each of these $k$ paths contains at most $\ell$ vertices. So $k + t - 1 \leqslant k\ell$ and thus $t \leqslant k(\ell - 1) + 1 \leqslant p(\ell - 1) + 1$.

Now let $t$ be the maximum integer such that $P$ is $t$-rotund. Then $P$ is not $(t + 1)$-rotund. Since $G$ is $\ell$-path free, by the analysis above, $s \leqslant 2s^* \leqslant 2[(t + 1)^p - 1] \leqslant 2[p(\ell - 1) + 2]^p - 2 \leqslant 2(\ell^2 - \ell + 2)^\ell - 2$.
\qed

A \emph{rooted hypergraph} is a hypergraph $G$ with a special designated subset $r(G)$ of $V(G)$. Let $\mathbf{Q}$ be a set with a quasi-ordering $\leqslant_{\mathbf{Q}}$. A \emph{$\mathbf{Q}$-labeled rooted hypergraph} is a rooted hypergraph $G$ with a mapping $\sigma: E(G) \mapsto \mathbf{Q}$.
\begin{lem}\label{lem_nboundedbqo}
Let $\mathbf{Q}$ be a better-quasi-ordered set, and $\mathcal{G}$ be a sequence of $\mathbf{Q}$-labeled rooted hypergraphs (respectively, of bounded multiplicity) with vertex sets the subsets of $[p]$, where $p \geqslant 1$. For $G_1, G_2 \in \mathcal{G}$, denote by $G_1 \subseteq G_2$ (respectively, $G_1 \leqslant G_2$) that $r(G_1) = r(G_2)$, and there is an isomorphism $\varphi$ from $G_1$ to a (respectively, an induced) subgraph  of $G_2$ such that for $i \in V(G_1)$ and $e \in E(G_1)$, $\varphi(i) = i$ and $\sigma(e) \leqslant_{\mathbf{Q}} \sigma(\varphi(e))$. Then $\mathcal{G}$ is better-quasi-ordered by $\subseteq$ (respectively, $\leqslant$).
\end{lem}
\pf
There are $\sum_{i = 0}^p\binom{p}{i}2^i = 3^p$ choices for vertex sets and roots. By Lemma \ref{lem_coverbqo}, it is safe to assume that all $G \in \mathcal{G}$ have the same vertex set, say $[p]$, and the same root. Then each $G$ is a sequence of length $2^p - 1$, indexed by the nonempty subsets of $[p]$, on the set of sequences of $\mathbf{Q}$. By Lemma \ref{lem_coverbqo} and \ref{lem_subseqbqo}, $\mathcal{G}$ is better-quasi-ordered by $\subseteq$. Now let $\mu$ be an upper bound of the multiplicity. There are $(\mu + 1)^{2^p - 1}$ unequal hypergraphs of vertex set $[p]$. By Lemma \ref{lem_coverbqo}, we can assume that all these rooted hypergraphs are equal. In this situation, each $G \in \mathcal{G}$ is a sequence of length $2^p - 1$ on the set of sequences of $\mathbf{Q}$ of length $\mu$. By Lemma \ref{lem_coverbqo} and \ref{lem_subseqbqo}, $\mathcal{G}$ is better-quasi-ordered under $\leqslant$.
\qed

Now we show that, for a better-quasi-ordered set $\mathbf{Q}$, the $\mathbf{Q}$-labeled hypergraphs of bounded (respectively, multiplicity,) tree-width and tree-diameter are better-quasi-ordered by (respectively, induced) subgraph relation.
\begin{lem}\label{lem_Qtwtdbqo}
Let $p, s \geqslant 0$ be integers, $\mathbf{Q}$ be a better-quasi-ordered set, $\mathcal{G}$ be the set of $\mathbf{G} := (G, T, \mathcal{V}, r, V_G)$, where $G$ is a $\mathbf{Q}$-labeled  hypergraph (respectively, of bounded multiplicity) with a tree-decomposition $(T, \mathcal{V})$ of width at most $p - 1$, $T$ is a rooted tree of root $r$ and height at most $s$, and $V_G \subseteq V_r$. Let $\lambda: V(G) \mapsto [p]$ be a colouring such that for each $v \in V(T)$, every pair of different vertices of $V_v$ are assigned different colours. For $\mathbf{X}, \mathbf{Y} \in \mathcal{G}$ (respectively, $\mathcal{G}^*$), denote by $\mathbf{X} \subseteq \mathbf{Y}$ (respectively, $\mathbf{X} \leqslant \mathbf{Y}$) that there exists an isomorphism $\varphi$ from $X$ to a subgraph (respectively, an induced subgraph) of $\mathbf{Y}$ such that $\varphi(V_X) = V_Y$, and that for each $x \in V(X)$ and $e \in E(X)$, $\lambda(x) = \lambda(\varphi(x))$ and $\sigma(e) \leqslant \sigma(\varphi(e))$. Then $\mathcal{G}$ is better-quasi-ordered by $\subseteq$ (respectively, $\leqslant$).
\end{lem}
\pf
Let $\mathcal{G}_s$ be the set of $\mathbf{G} \in \mathcal{G}$ of which the height of $T$ is exactly $s$. By Lemma \ref{lem_coverbqo}, it is enough to prove the lemma for $\mathcal{G}_s$. The case of $s = 0$ is ensured by Lemma \ref{lem_nboundedbqo}. Inductively assume it holds for some $s - 1 \geqslant 0$. By Lemma \ref{lem_subseqbqo}, the powerset $\mathcal{M}_{s - 1}$ of $\mathcal{G}_{s - 1}$ is better-quasi-ordered.
For each $u \in N_T(r)$, let $T_u$ be the component of $T - r$ containing $u$, $G_{T_u}$ be the subgraph of $G$ induced by the vertex set $\bigcup_{w \in T_u}V_w$, $\mathcal{V}_{T_u} := \{V_w | w \in V(T_u)\}$, and $\mathbf{G}_{T_u} := (G_{T_u}, T_u, \mathcal{V}_{T_u}, u, V_r \cap V_u) \in \mathcal{G}_{s - 1}$. Let $G_r$ be the subgraph of $G$ induced by $V_r$, and $\mathbf{G}_r := (G_r, r, V_r, r, V_G)$. Clearly, $\mathbf{G} \mapsto \mathbf{G}_r \times \{\mathbf{G}_{T_u}| u \in N_T(r)\}$ is an order-preserving bijection from $\mathcal{G}_s$ to $\mathcal{G}_0 \times \mathcal{M}_{s - 1}$. By Lemma \ref{lem_coverbqo}, $\mathcal{G}_s$ is better-quasi-ordered since $\mathcal{G}_0$ and $\mathcal{M}_{s - 1}$ are better-quasi-ordered.
\qed

We are ready to prove the rest of our main results.

\noindent {\bf Proof of Lemma \ref{lem_lpathfreebqo}.} $(\Leftarrow)$ follows from Lemma \ref{lem_tdibounded} and \ref{lem_Qtwtdbqo}. $(\Rightarrow)$ Suppose for a contradiction that $H$ is not a union of paths. Then $H$ contains a cycle or a vertex of degree at least $3$. Let $n := n(H)$, and $O_i$ be the cycle of $n + i$ vertices. Then $O_1, O_2, \ldots, $ is a bad sequence with respect to subgraph or induced subgraph relation.\qed

\noindent {\bf Proof of Theorem \ref{thm_pathrootsBound}.} The upper bounds for the order, size and total number of minimal $\ell$-path roots of a finite graph $H$ are given in Lemma \ref{lem_pathrootosnbound}. Let $n := n(H)$. For each $G \in \mathbb{Q}_\ell[H]$, suppose for a contradiction that $G$ contains a path of length $\ell n$. Then $H$ contains an $n$-path with $n + 1$ vertices, a contradiction. The rest of the theorem follows from Lemma \ref{lem_lpathfreebqo}.\qed